\documentclass[12pt]{article}
\usepackage{amsmath,amssymb,amsfonts,amsthm}
\newcounter{item}[section]
\newcounter{kirshr}
\newcounter{kirsha}
\newcounter{kirshb}
\newenvironment{enumroman}{\setcounter{kirshr}{1}
\begin{list}{(\roman{kirshr})}{\usecounter{kirshr}} }{\end{list}}
\newenvironment{enumarab}{\setcounter{kirshb}{1}
\begin{list}{(\arabic{kirshb})}{\usecounter{kirshb}} }{\end{list}}
\newenvironment{athm}[1]{\vskip3mm\par\noindent
{\bf #1 }. \slshape }
{\upshape\par\vskip10pt minus3pt}
\newtheorem{theorem}{Theorem}[section]

\newtheorem{lemma}[theorem]{Lemma}
\newtheorem{corollary}[theorem]{Corollary}
\newenvironment{demo}[1]{\noindent{\bf #1.}\upshape\mdseries}
{\nopagebreak{\hfill\rule{2mm}{2mm}\nopagebreak}\par\normalfont}
\theoremstyle{definition}

\newtheorem{example}[theorem]{Example}
\newtheorem{definition}[theorem]{Definition}

\def\C{{\mathfrak{C}}}
\def\Fm{{\mathfrak{Fm}}}

\def\Nr{{\mathfrak{Nr}}}
\def\Fr{{\mathfrak{Fr}}}
\def\Sg{{\mathfrak{Sg}}}
\def\Fm{{\mathfrak{Fm}}}
\def\A{{\mathfrak{A}}}
\def\B{{\mathfrak{B}}}
\def\C{{\mathfrak{C}}}
\def\D{{\mathfrak{D}}}
\def\M{{\mathfrak{M}}}

\def\F{{\mathfrak{F}}}

\def\Rd{{\ Rd}}
\def\(R)RA{{\bf (R)RA}}

\def\B{{\sf B}}

\def\g{{\sf g}}

\def\Nr{{\mathfrak{Nr}}}

\def\Nr{{\mathfrak{Nr}}}

\def\A{{\mathfrak{A}}}
\def\B{{\mathfrak{B}}}
\def\C{{\mathfrak{C}}}
\def\D{{\mathfrak{D}}}

\def\A{{\mathfrak{A}}}
\def\B{{\mathfrak{B}}}
\def\C{{\mathfrak{C}}}
\def\D{{\mathfrak{D}}}

\def\L{{\mathfrak{L}}}
\def\Rd{{\mathfrak{Rd}}}

\def\L{{\mathfrak{L}}}

\title{Interpolation in many valued logics using algebraic logic}
\author{Tarek Sayed Ahmed \\
Department of Mathematics, Faculty of Science,\\ 
Cairo University, Giza, Egypt.\\
 E mail: rutahmed@gmail.com\\ 
Fax: 02-33381411}
\begin{document}
\maketitle

\begin{abstract}
\noindent  We prove several interpolation theorems for many valued infinitary logic with quantifiers, 
by studying expansions of $MV$ algebras in the spirit of polyadic and cylindric algebras.

\end{abstract}

$MV$ algebras were introduced by Chang in 1958 \cite{Chang} to provide an algebraic reflection of the completeness theorem of 
the Lukasiewicz infinite valued propositional logic. In recent years the range of applications of $MV$ algebras
has been enormously extended with profound interaction with other topics, ranging from  
lattice ordered abelian groups, $C^*$ algebras, to fuzzy logic.
In this paper we study $MV$ algebras in connection to fuzzy (many valued) logic. 
We prove five interpolation theorems for many valued logic using the machinery 
of algebraic logic, four interpoltation theorems for fuzzy logic in the narrow sense and one for rational Pavelka logic where truth values are incorporated in the syntax.
One of the earliest papers (if not the first) that deals with an application of polyadic algebras to $MV$ algebras is \cite{Sh}.
Another application of polyadic algebras to $MV$ algebras and Pavelka Rational logic is \cite{Pav}. 
Here our proven interpolation theorems substantially generalize the
representability results proved in the formentioned papers.

An $MV$ algebra, has a dual behaviour; it  can be viewed,  
in one of its  facets, as  a `non-idempotent' generalization of a Boolean algebra possesing a strong lattice structure. 
The lack of idempotency enables $MV$ algebras  to be compared to monodial structures like monoids and abelian groups. 
Indeed, the category of $MV$ algebras has been shown to be equivalent to the category of $l$ groups. At the same time the lattice structure 
of Boolean algebras can be recovered inside $MV$ algebras, by an appropriate term definability of primitive connectives.  In this 
respect, they  have a strong lattice structure 
(distributive and bounded), which make the techniques of lattice theory readily applicable to their study.
As shown in this paper, in certain contexts when we replace the notion of a Boolean algebra with an $MV$ algebra, 
the results survive such a replacement with
some non-trivial modifications, and this can be accomplished in a somewhat  unexpected manner. 

Boolean algebras work as the equivalent semantics of classical propositional logic. To sudy classical first order logic, 
Tarski \cite{HMT1}, \cite{HMT2} introduced cylindric algebras, while Halmos \cite{Halmos} introduced polyadic algebras. 
Both of those can be viewed as Boolean algebras with extra operations that reflect algebraically 
existential quantifiers. 

Boolean algebras also have a neat and intuitive depiction, modulo isomorphisms; 
any Boolean algebra is an algebra of subsets of some set endowed 
with the concrete set theoretic operations of union, intersection and complements. Such a connection, a typical duality theorem, 
is today well understood. These nice properties mentioned above is formalized through the topology of 
Stone spaces that allows to select the right objects in the full power set of some set, the underlying set of the 
associated topological space. The representation theory of cylindric algebras, on the other hand, proves much more involved, 
 and lacks such a strong well understood duality theorem like that of 
Boolean algebras. However, there is an extension of Stone duality to cylindric algebras, due to Comer \cite{C}, where he establishes a 
dual equivalence between cylindric algebras 
and certain categories of sheaves; but such a duality does not go deeper into the analysis of representability.
There is a version of concrete  (representable) algebras for cylindric algebras, with extra operations interpreted as projections,
but this does not coincide with the abstract class of cylindric algebras. This is in sharp contrast to 
Boolean algebras. It is not the case 
that every cylindric algebra is representable in a concrete manner with the operations being set theoretic operations on relations. Not only that, 
but in fact the class of representable algebras need an infinite axiomatization in first order logic, and for any such axiomatization, 
there is an inevitable degree of complexity \cite{Andreka}. 
On the other hand, polyadic algebras enjoy a strong representation theorem; every polyadic algebra is representable \cite{DM}.
Here we apply the theory of polyadic algebras to $MV$ algebras. The idea is to study transformation systems based on such algebras.

A polyadic algebra is typically an instance of a transformation system. A transformation system can be defined to be a quadruple of 
the form $(\A, I, G, S)$ where $\A$ is an algebra (in case of polyadic algebras it is a Boolean algebra), $I$ is a non empty set (we will only be concerned with infinite sets)
$G$ is a subsemigroup of $(^II,\circ)$ and $S$ is a homomorphism from $G$ to the semigroup of endomorphisms of $\A$, denoted by $End(\A)$. 
Elements of $G$ are called transformations. 
For a transformation $\tau$, $S({\tau})\in End(\A)$ is  called a substitution. 
Polyadic algebras arise when $\A$ is a Boolean algebra endowed with quantifiers and $G={}^II$. 
There is an extensive literature for polyadic algebras dating back to the fifties and sixties of the last century, 
\cite{Halmos}, \cite{J}, \cite{D}, \cite{DM}, \cite{ AUamal}, \cite{S}. 
Introduced by Halmos in the fifties of the last century, the theory of polyadic algebras is now picking up again; indeed it's
regaining momentum with pleasing progress and a plathora of results, to mention a few references 
in this connection, the reader is referred to  \cite{NS}, \cite{MLQ}, \cite{Fer1}, \cite{Fer2}, \cite{Fer3}, \cite{Fer4}, \cite{ANS}, \cite{trans}.
In recent times reducts of polyadic algebras were studied \cite{S}, \cite{AUamal}; these reducts 
are obtained by restricting quantifiers to involve only quantification on finitely many 
variables and to study 
(proper) subsemigroups of $^II$ (endowed with the binary operation of composition of maps.) 
The two extremes are the semigroup of finite transformations 
(a finite transformation is one that moves only finitely many points)
and all of $^II$ but there are infinitely many semigroups in between.

In this paper we study reducts of polyadic algebras by allowing (proper) subsemigroups of $^II$, 
but we also weaken the Boolean structure  to be an $MV$ algebra. Thus we are in the realm of many valued  quantifier infinitary logics.
Many-valued logics are non-classical logics. The most two basic (semantical) assumptions of classical logic propositional as well as first order 
are the principles of bivalence and of compositionality. The principle of bivalence is the assumption that each sentence is either true or false, 
i.e has exactly one truth value. Many valued logics differ from classical ones  by the fundamental fact that it does  
not restrict the number of truth values to only two: they allow for a larger set (possibly infinite)
of truth degrees.
However, like classical logic they accept the 
principle of compostionality (or truth-functionality), namely, that the truth of a compound sentence is determined by 
the truth values of its component sentences 
(and so remains unaffected when one of its component sentences is replaced by another sentence with the same truth value).

The formalized languages for systems of many-valued logic follow the two standard patterns for propositional and predicate logic, respectively:
there are propositional variables together with connectives and (possibly also) truth degree constants in the case of propositional languages, 
there are object variables together with predicate symbols, possibly also object constants and function symbols, 
as well as quantifiers, connectives, and (possibly also) truth degree constants in the case of first-order languages.
We shall deal with both cases, infinitary extensions of Lukasiewicz  predicate logic and the Pavelka  Rational logics, using polyadic algebras 
where transformations are restricted to a 
semigroup $G$.

We shall study the cases when $G$ consists of all finite transformations on a set $I$, when $G$ 
is a proper subsemigroup of $^II$ satisfying certain properties but essentially containing infinitary substitutions (this involves infinitely many cases), 
and when $G$ is the semigroup of all transformations. In all these 
three  cases the scope of quantifiers are finite, so in this respect our algebras also 
resemble cylindric algebras. However, the fourth case we study is that when $G={}^II$ and the scope of quantifiers is infinite. 
So in the latter case we are in the polyadic paradigm.

It is folklore in algebraic logic that cylindric algebras and polyadic algebras belong to different paradigms.
For example the class of representable cylindric algebras admits a recursive axiomatization, 
but all axiomatizations of polyadic algebras are extremely complex from the recursive point of view; 
they are  not recursively enumerable.
The interaction between the theories of cylindric algebras and polyadic algebras is extensively 
studied in algebraic logic, see e.g
\cite{ANS}, with differences and similarities illuminating both theories.
In fact the study of $G$ Boolean polyadic algebras by Sain in her pioneering paper \cite{S} is an outcome of such research; 
it's a typical situation in which the positive properties of both theories amalgamate.
Boolean polyadic algebras, when $G$ is the set of finite transformations of $I$ or $G={}^II$ are old, in this respect the reader is referred to 
\cite{Halmos}, \cite{D},
\cite{DM}. In the former case such algebras are known as quasipolyadic 
algebras, which are substantially different than polyadic algebras (in the infinite dimensional case), as is well known, quasipolyadic algebras
belong to the cylindric paradigm.
However, studying reducts of polyadic algebras by allowing only 
those substitutions coming from a subsemigroup of $^II$ is relatively recent starting at the turn of the last century
\cite{S}.

Such  algebras, with Boolean reducts (of which we study their $MV$ reducts), also provide a possible  solution to a central 
problem in algebraic logic, better known as the finitizability problem,  which asks for a simple (hopefully) finite axiomatization 
for the class of representable algebras.  The class of representable algebras is given by specifying the universes of the algebras in the 
class, as sets of certain sets
endowed with set theoretic concrete operations; thus representable algebras are completely determined once one specifies their universes.
The finitizability problem - the attempt to capture the essence of such set-theoretic constructions in a thoroughly finitary matter - proved to be a difficult 
task, and has been 
discussed at length  in the literature \cite{Bulletin}. We shall show that the $MV$ polyadic reducts of such algebras, 
still form a finitely axiomatizable variety, consisting solely of representable algebras, that further enjoys a strong form of amalgamation
(known as the superamalgamation property), a result that could be of interest to abstract algebraic logic, as well.

Being rather a family of problems, the finitizability problem  
has several scattered reincarnations and in some sense is still open \cite{Simon}. 
The finitizability problem has philosophical implications concerning reasoning about reasoning, and can be likened to Hilbert's programe
of proving mathematics consistent by concrete finitistic methods.

Heyting polyadic algebras, as transformations systems, were studied by Monk \cite{Monk}, Georgescu \cite{G} and the present author \cite{Hung}. 
We continue this trend by studing expansions of $MV$ polyadic algebras as transformation systems.

We prove algebraically, using the techniques of Halmos and Tarski in algebraic logic,  four interpolation theorem  which show that predicate 
Lukasiewicz calculas together with some of its infinitary extensions, have the Craig interpolation property. 
For many valued predicate logics, the  main types of logical calculi are Hilbert style calculi, Genzen type sequent calculi and Tableau calculi.
From our result, it readily follows that  the extensions we study are also complete, relative to a Hilbert style calculas.
Interpolation property for various logics both classical and non-classical have been studied extensively in the literature of algebraic logic, 
\cite{AUU}, \cite{IGPL}, \cite{MLQ}, \cite{Hung}, \cite{D}, \cite{J}, \cite {Shelah}, \cite{Mak}.

A historic comment is in order. Formalized many valued logic can be traced back  to the work of Lukasiewicz in 1920 
and the independent work of Post in 1921, when three valued is studied. 
Heyting, a few years later, introduced a three valued propositional calculus related to intuitionistic logic.
G\'odel proposed an infinite hierarchy of finitely-valued systems; 
his goal was to show that intuitionistic logic is not a many valued logic.
In the last few decades many valued logics have acquired tremendous interest; 
in 1965 Zadeh had published his landmark paper on fuzzy sets and the trend of fuzzy Logic started. 
Today the various approaches to many valued logics are aspiring to provide fuzzy logic the theoretical rigorous 
foundations that were lacking  for a long time.

Throughout the paper, we follow more or less standard notation. We distinguish notationally between an algebra $\A$ and its universe $A$. 
When we write Gothic letters for algebras $\A,\B, \ldots$, it is to be tacitly assumed that the corresponding Roman letters $A,B,\ldots$ denote 
their universes. Otherwise, unfamiliar notation will be introduced at their first occurrance in the text. 

\section{Prelimenaries and the main results in logical form}

\begin{definition} An $MV$ algebra is an algebra
$$\A=(A, \oplus, \odot, \neg, 0,1)$$
where $\oplus$, $\odot$ are binary operations, $\neg$ is a unary operation and $0,1\in A$, such that the following identities hold: 
\begin{enumerate}
\item $a\oplus b=b\oplus a,\ \ \  a\odot b=b\odot a.$
\item $a\oplus (b \oplus c)=(a\oplus b)\oplus c$,\  \  \ $a\odot (b \odot c)=(a\odot b)\odot c.$

\item $a\oplus 0=a$ ,\ \ \ $a\odot 1=a.$
\item $a\oplus 1=1$, \ \ \ $a\odot 0=a.$
\item  $a\oplus \neg a=1$,\ \ \  $a\odot \neg a=0.$
\item $\neg (a\oplus b)=\neg a\odot \neg b,$\ \ \ $\neg (a\odot b)=\neg a\oplus \neg b.$
\item $a=\neg \neg a$\ \ \ $\neg 0=1.$
\item $\neg(\neg a\oplus b)\oplus b=\neg(\neg b\oplus a)\oplus a.$
\end{enumerate}
\end{definition}
$MV$ algebras form a variety that is a subvariety of the variety of  $BL$ algebras intoduced by Hajek, 
in fact $MV$ algebras coincide with those $BL$ algebras satisfying double negation law, 
namely that $\neg\neg x=x$, and contains all  Boolean algebras.
\begin{example} A simple numerical example is $A=[0,1]$ with operations $x\oplus y=min(x+y, 1)$, $x\odot y=max(x+y-1, 0)$,  and $\neg x=1-x$. 
In mathematical fuzzy logic, this $MV$-algebra is called the standard $MV$ algebra, 
as it forms the standard real-valued semantics of Lukasiewicz logic.
\end{example}
$MV$ algebras aso arise from the study of continous $t$ norms. 
\begin{definition}A $t$ norm is a binary operation $*$ on $[0,1]$, i.e $(t:[0,1]^2\to [0,1]$) such that
\begin{enumroman}
\item  $*$ is commutative and associative,
that is for all $x,y,z\in [0,1]$,
$$x*y=y*x,$$
$$(x*y)*z=x*(y*z).$$
\item $*$ is non decreasing in both arguments, that is
$$x_1\leq x_2\implies x_1*y\leq x_2*y,$$
$$y_1\leq y_2\implies x*y_1\leq x*y_2.$$
\item $1*x=x$ and $0*x=0$ for all $x\in [0,1].$

\end{enumroman}
\end{definition}
The following are the most important (known) examples of continuous $t$ norms.

\begin{enumroman}
\item Lukasiewicz $t$ norm: $x*y=max(0,x+y-1),$
\item Godel $t$ norm $x*y=min(x,y),$
\item Product $t$ norm $x*y=x.y$.
\end{enumroman}
We have the following known result \cite{H} lemma 2.1.6
 
\begin{theorem} Let $*$ be a continuous $t$ norm. 
Then there is a unique binary operation $x\to y$ satisfying for all $x,y,z\in [0,1]$, the condition $(x*z)\leq y$ iff $z\leq (x\to y)$, namely 
$x\to y=max\{z: x*z\leq y\}.$
\end{theorem}
The operation $x\to y$ is called the residuam of the $t$ norm. The residuam $\to$ 
defines its corresponding unary operation of precomplement 
$\neg x=(x\to 0)$.
Abstracting away from $t$ norms, we get

\begin{definition} A residuated lattice is an algebra
$$(L,\cup,\cap, *, \to 0,1)$$ 
with four binary operations and two constants such that
\begin{enumroman}
\item $(L,\cup,\cap, 0,1)$ is a lattice with largest element $1$ and the least element $0$ (with respect to the lattice ordering defined the usual way: 
$a\leq b$ iff $a\cap b=a$).
\item $(L,*,1)$ is a commutative semigroup with largest element $1$, that is $*$ is commutative, associative, $1*x=x$ for all $x$.
\item Letting $\leq$ denote the usual lattice ordering, we have $*$ and $\to $ form an adjoint pair, i.e for all $x,y,z$
$$z\leq (x\to y)\Longleftrightarrow x*z\leq y.$$
\end{enumroman}
\end{definition}
A result of Hajek, is that an $MV$ algebra is  a prelinear commutative bounded integral residuated lattice 
satisfying the additional identity $x\cup y=(x\to y)\to y.$ In case of an $MV$ algebra, $*$ is the so-called strong conjunction which we 
denote here following standard notation in the 
literature by $\odot$. $\cap$ is called weak conjunction. The other operations are defined by 
$\neg a=a\to 0$ and $a\oplus b=\neg(\neg a\odot \neg b).$ The operation $\cup$ is called weak disjunction, while $\oplus$ is called 
strong disjunction. The presence of weak and strong conjunction is a common feature of substructural logics without the rule of contraction, to which Lukasiewicz 
logic belongs.

We now turn to describing some metalogical notions, culminating in formulating our main results in logical form. However, throughout the paper, 
our investigations will be purely algebraic,
using the well develped machinery of algebraic logic. 
There are two kinds of semantics for systems of many-valued logic.
Standard logical matrices and algebraic semantics. We shall only encounter  algebraic semantics. 
From a philosophical, especially epistemological point of view the semantic aspect of logic is more basic than the syntactic one,  
because it is mainly the semantic core which determines the choice of suitable 
syntactic versions of the corresponding system of logic.

Informally, a language is a triple $\Lambda=(V, P, G)$ where $V$ is a set providing an infinite supply of variables, 
$P$ is a another set of predicates disjoint from $V,$
and $G$ is a semigroup of transformations on $V$. 
There is no restriction on the arity of $p\in P$, that is the arity may be infinite. 
Formulas are defined recursively the usual way.  Atomic formulas are of the form $p\bar{v}$, the length of $\bar{v}$ is equal to the arity of $p$.
If $\phi, \psi$ are formulas, $W\subseteq V$ and $\tau\in G$, 
then $\phi\oplus\psi$, $\phi\odot \psi$, $\neg\phi$, $\exists W\phi,$ $\forall W\phi$, and ${\sf S}({\tau})\phi$ are formulas. 
Since we allow infinitary predicates, the scope of quantification can be infinite. 

A structure for a predicate language is $\M=(M, p_M)_{p\in P}$ 
where $M\neq \emptyset$, for each predicate
$p$ of arity $n$, $n$ an ordinal (could be infinite), $p_M$ is an $n$-ary $[0,1]$ 
fuzzy relation on $M,$ that is $p_M:{}^nM\to [0,1]$.
For each formula $\phi$ the truth value $||\phi||_{\M,s}$ of $\phi$ in $\M$ is determined by the  
evaluation $s$ of free variables the usual 
Tarskian way.

In more detail, an $\M$ evaluation is a map from $V$ to $M$. For two evaluations $s$ and $s'$ and $\Gamma\subseteq V$, we write 
$s\equiv_{\Gamma} s'$ iff $s(v)=s'(v)$ for all $v\notin \Gamma$.
The value of a variable given by $\M,s$ is defined by $||v||_{\M,s}=s(v)$. For a formula $\phi$ and a transformation $\tau\in {}^VV$, we write
${\sf S}(\tau)\phi$ for the formula obtained by the simultaneous substitution of the variable $v_{\tau(i)}$ for $v_i$
such that the substitution is free. In presence of infinitary quantification, there could be a risk of collision of variables, 
but we will indicate below how such a possibly 
infinitary operation can be always executed. 

Now we define the truth value $||\phi||_{\M,s}$:
$$||p_M(v_1\ldots v_i\ldots )||_{\M,s}=p_M(s(v_1)\ldots s(v_i)\ldots ),$$
$$||\phi\oplus \psi||_{\M,s}=||\phi||_{\M,s}\oplus ||\psi||_{\M,s},$$
$$||\phi\odot\psi||_{\M,s}=||\phi||_{\M,s}\odot ||\psi||_{\M,s},$$
$$||\neg \phi||_{\M,s}=\neg ||\phi||_{\M,s},$$
$$||{\sf S}(\tau)\phi)||_{\M,s}=||\phi||_{\M, s\circ \tau},$$
$$||(\exists W)\phi||_{\M,s}=\bigvee\{||\phi||_{\M,s'}: s'\equiv_W s\}.$$
$$||(\forall W)\phi||_{\M,s}=\bigwedge\{||\phi||_{\M,s}: s'\equiv_W s\}.$$
A structure is safe if the suprema and infima in the last two clauses always exist. Thus all structures considered are safe 
(since $[0,1]$ with the usual order is a complete lattice).
Let $\M=(M, p_M)_{p\in P}$ be a structure and $s:V\to M$.
A formula $\phi$ is satisfiable under $s$ if $||\phi||_{\M,s}=1$. $\phi$ is valid in $\M$ and $\M$ is a model of $\phi$ if $||\phi||_{\M,s}=1$ for 
all $s\in {}^VM$.
We write $\M\models \phi$ if $\phi$ is valid in $\M$. For a set of formulas $\Gamma$, we 
write $\Gamma\models \phi$ if for every structure $\M$ whenever every formula in $\Gamma$ is valid in $\M$, then $\phi$ 
itself is valid in $\M$. We write $\models \phi$ for $\emptyset\models \phi$.

Since we dealing with infinitary language we need to formalize such languages rigorously in set theory.
We now define a calculus (in a usual underlying set theory $ZFC$,  say) that we prove to be complete with respect to usual semantics; 
this will follow from our stronger proven result that such logics 
enjoy the interpolation property. 
Let $V$ and $P$ be disjoint sets of symbols, such that $V$ is infinite, $\rho$ a function 
with domain $P$ and whose range is a set of ordinals. Let $\mathfrak{m}$ be  a cardinal, such that 
$\mathfrak{m}\leq |V|^+$. (For a cardinal $\mathfrak{n}, \mathfrak{n}^+$ denotes its successor). 
Assume in addition that for each $p\in P,$ $|\rho(p)|\leq |V|$.
For a set $I$, let $S_{\mathfrak{m}}I=\{X: X\subseteq I: |X|<\mathfrak{m}\}$.
Fix $G\subseteq {}^VV$ and $T\subseteq S_{\mathfrak{m}}(V)$.
We define (formally) a logic $\mathfrak{L}_{V, \rho, \mathfrak{m},G, T}$, 
based on $G$, or simply $\mathfrak{L}_{G,T}$\footnote{
We could omit some of the subscripts and keep others that are relevant to the context, for example, we may 
write $\mathfrak{L}_{\mu,G,T}$, when the other superscripts are immaterial in the given context.}  in the following way.
The symbols of $L$ consist of: 
\begin{enumarab}
\item the (strong) disjunction $\oplus$, (strong) conjunction $\odot $, implication symbol $\to$, and negation symbol $\neg$,
\item two symbols $\bot$, $\top$ standing for falsity and truth,
\item universal quantification symbol $\forall$,
\item existential quantification symbol $\exists$,
\item the individual variables $v\in V$ and predicates $p\in P,$
\end{enumarab}
We note that $\forall$ and $\to$ are redundant, but their introduction from the start, makes life easior.
We assume that $\oplus, \odot,\to, \neg, \forall, \exists$ are not members of $V$ nor $P$.
An atomic formula is an ordered pair $(p,x)$ 
where $p\in P$ and $x\in {}^{\rho(p)}V.$ Formulas are defined the usual way by recursion:
$\phi$ is a formula in $\mathfrak{L}_{G,T}$ 
if there exists a finite sequence $\phi_0,\phi_1,\ldots, \phi_n$ called a formation of $\phi$ in $\mathfrak{L}_{G,T}$, 
such that $\phi_n=\phi$ and for each $m\leq n$, at least one of the following conditions hold:
\begin{enumarab}
\item $\phi_m$ is an atomic formula in $\mathfrak{L}_{G,T},$
\item $\phi_m$ is either $\top$ or $\bot,$
\item for some $k,l<m$, $\phi_m$ is one of the ordered triplet $(\phi_k,\oplus, \phi_l)$,  $(\phi_k,\odot, \phi_l)$ or  
$(\phi_k,\to, \phi_l),$
\item for some $k<m$, $\phi_m$ is the ordered pair $(\neg, \phi_l),$
\item for some $k<m$, and $W\in T$, $\phi_m$ is the ordered triplet $(\forall, W, \phi_k),$
\item for some $k<m$, and $W\in T$, $\phi_m$ is the ordered triplet $(\exists, W, \phi_k).$
\end{enumarab}
We often omit commas in ordered pairs or triplets or quadruples, so we may write
$(\exists W\phi)$\footnote{This is a notation, it should not be mistaken for brackets, which do not exist in our vocabularly.} for $(\exists, W, \phi)$.
The set $V_f(\phi)$ of free variables and the set $V_b(\phi)$ of bound variables in a formula $\phi$ are defined recursively the usual way.
That is 
\begin{enumarab}
\item If $\phi$ is an atomic formula $(px)$, then $V_f(\phi)$ is the range of $x$.
\item If $\phi$ is $(\psi\oplus  \theta)$ or $(\psi\odot \theta)$ or $(\psi\to \theta)$, then $V_f(\phi)=V_f(\psi)\cup V_f(\theta).$
\item If $\phi$ is $(\neg \psi)$ then $V_f(\phi)=V_f(\psi)$.
\item
If $\phi=(\forall W \psi)$ or $(\exists W\psi)$,  then $V_f(\psi)=V_f(\phi)\sim W$.
Now for the bound variables $V_b(\phi)$:
\item If $\phi$ is an atomic formula $(px)$, then $V_b(\phi)=0.$ 
\item If $\phi$ is $(\psi\oplus  \theta)$ or $(\psi\odot  \theta)$ or $(\psi\to \theta)$,  then $V_b(\phi)=V_b(\psi)\cup V_b(\theta).$
\item  if $\phi$ is $\neg \psi$, then $V_b(\phi)=V_b(\psi)$.
\item If $\phi=(\forall W \psi)$, then $V_b(\psi)=V_b(\phi)\cup W.$
\item If $\phi=(\exists W \psi)$, then $V_b(\psi)=V_b(\phi)\cup W.$
\end{enumarab}
In the two cases of $\phi=\top$ or $\bot$, then $V_f(\phi)=V_b(\phi)=0$.
Note that the variables occuring in a formula $\phi$, denoted by $V(\phi),$  is equal to $V_f(\phi)\cup V_b(\phi)$ which could well be infinite.
For $\tau\in G$ and $\phi$ a formula, ${\sf S}(\tau)\phi$ (the result of substituting each variable $v$ in $\phi$ by $\tau(v)$)
is defined recursively and so is ${\sf S}_f(\tau)\phi$ (the result of substituting each free variable $v$ by $\tau(v)$).\footnote{Notice that ${\sf S}(\tau)$ 
and for that matter $S_f(\tau)$ is not part of the vocabularly of our languages. We do not have transformations in $G$ as symbols.}
\begin{enumarab}
\item If $\phi$ is atomic formula $(px)$ then ${\sf S}(\tau)\phi=(p,\tau\circ x).$
\item If $\phi$ is $(\psi\oplus \theta)$ then  ${\sf S}(\tau)\phi=({\sf S}(\tau)\psi\oplus {\sf S}(\tau)\theta).$ The same for other propositional connectives.
\item If $\phi=(\forall W\phi)$ then ${\sf S}(\tau)\phi=(\forall\tau(W){\sf S}(\tau)\phi).$
\item If $\phi=(\exists  W\phi)$ then ${\sf S}(\tau)\phi=(\exists\tau(W){\sf S}(\tau)\phi).$

\end{enumarab}
To deal with free substitutions, we introduce a piece of notation that proves helpful. For any function $f\in {}^XY$ and any set $Z$, we let 
$$f|Z=\{(x, f(x)): x\in X\cap Z\}\cup \{(z,z)|z\in Z\sim X\}.$$
Then $f|Z$ always has domain $Z$ and $0|Z$ is the identity function on $Z$.

Now for free subtitutions  the first two  clauses are the same, but if
$\phi=(\forall W \psi)$, then ${\sf S}_f({\tau})\phi=(\forall W {\sf S}_f(\sigma)\psi)$ 
and if $\phi=(\exists  W \psi)$ then ${\sf S}_f({\tau})\phi=(\exists v {\sf S}_f(\sigma)\psi)$ 
where $\sigma=\tau|(V\sim W)|V$.

If $\tau\in \bigcup\{^WV: W\subseteq V\}$, and $\phi$ is a formula, let ${\sf S}(\tau)\phi={\sf S}(\tau|V)\phi$ and ${\sf S}_f(\tau)\phi=S_f(\tau|V)\phi$.
Now we specify the axioms and the rules of inference.
The axioms are:
\begin{enumarab}
\item Axioms for propositional $MV$ logic using $\odot$, $\oplus$ and $\neg$.
\item $(((\phi\to \psi)\to (\neg \phi\oplus \psi))\odot  ((\neg \phi\oplus  \psi) \to (\phi\to \psi))).$
(Informally $(\phi\to \psi)$ is equivalent to $(\neg \phi\oplus \psi))$.
\item $((\forall W(\phi\to \psi)\to (\phi\to \forall W \psi)))$ where $W\in S_{\mathfrak{m}}(V\sim V_{f} \phi).$
\item $((\forall W(\phi\to \psi)\to (\exists W\phi\to \psi)))$ where $W\in S_{\mathfrak{m}}(V\sim V_{f}\phi).$
\item $(\forall W\phi\to {\sf S}_f(\tau)\phi)$, when  $\tau\in {}^W(V\sim V_{b}\phi).$
\item $({\sf S}_f(\tau)\phi\to (\exists W \phi))$, when  $\tau\in {}^W(V\sim V_{b}\phi).$
\end{enumarab}
The rules are:
\begin{enumarab}
\item From $\phi$, $(\phi\to \psi)$ infer $\psi.$ (Modus ponens.)
\item From $\phi$ infer $(\forall W\phi).$ (Rule of generalization.)
\item From ${\sf S}_f(\tau)\phi$ infer $\phi$ whenever $\tau\in {}^{V_f(\phi)} (V\sim V_b(\phi))$ and $\tau$ is one to one. (Free substitution.)
\item From $\phi$ infer ${\sf S}(\tau)\phi$ whenever $\tau\in {}^{V(\phi)}V$ is one to one (Substitution.)
\end{enumarab}
Proofs are defined the usual way, we write $\Gamma\vdash \phi$, if there is a proof of $\phi$ from $\Gamma$, that is there is a finite sequence 
$\phi_0,\phi_1\ldots \phi_n$ such that $\phi_n=\phi$ and for every $l<n$, either 
$\phi_l\in \Gamma$, or $\phi_l$ is an axiom, or $\phi_l$ follows from preceding formulas in the sequence by application of one of the rules of inference. 
\begin{definition} \begin{enumarab}
\item A logic $\mathfrak{L}_{G,T}$ is complete iff $\Gamma\models \phi$ implies $\Gamma\vdash \phi$.
\item A logic $\mathfrak{L}_{G,T}$ has the interpolation property if 
whenever $\models \phi\to \psi$, then there is a formula $\theta$ in the common vocabularly of $\phi$ and $\psi$ such that 
$\models \phi\to \theta$ and $\models \theta \to \psi.$
\end{enumarab}
\end{definition}
We consider four logics, one of which is countable. 
\begin{theorem} Let $V, \rho, G, T, \mathfrak{m}$ be as specified above. Let $\L_{V,\rho,G,T,\mathfrak{m}}$ be the corresponding logic. Then if
\begin{enumarab}
\item $V$ and $P$ are disjoint countable sets, $\mathfrak{m}=\omega$, $G$ is a 
rich semigroup (to be defined below) of ${}^VV$,  $T=S_{\omega}V$,
or
\item $V$ and $P$ are disjoint infinite sets , $\mathfrak{m}=|V|^+$, $G$ is the semigroup of finite transformations,
$T=S_{\omega}V$ and  $|V\sim \rho(p)|\geq \omega$,
or
\item $V$ and $P$ are disjoint infinite sets, $\mathfrak{m}=|V|^+$, $G$ is the semigroup of all transformations
 and $T=S_{\omega}V$, or

\item $V$ and $P$ are disjoint infinte sets, $\mathfrak{m}=|V|^+$, $G$ is the semigroup of all transformations
 and $T=S_{\mathfrak{m}}V$.

\end{enumarab}

Then $\mathfrak{L}_{V, \rho,G,T,\mathfrak{m}}$ is strongly complete and has the Craig interpolation property.
\end{theorem}
From $(2)$ when $\rho(p)$ is finite for every $p\in P$, then this gives an interpolation theorem for usual 
predicate Lukasiewicz logic.

\section{Algebraic Preliminaries}

For an algebra $\A$, $End(\A)$ denotes the set of endomorphisms of $\A$, i.e homomorphisms from $\A$ into itself.

\begin{definition} A transformation system is a quadruple $(\A, I, G, S)$ where $\A$ is an algebra, $I$ is a set, $G$ is a subsemigroup of $^II$ 
and $S$ is a homomorphism from
$G$ into $End(\A)$
\end{definition}
We shall deal with three cases of $G$. When $G$ is the semigroup of finite transformations, $G$ is a countable subset of $^II$ 
satisfying certain conditions but containing infinitary substitutions and $G={}^II.$
$\A$ will always be an $MV$ algebra.
If we want to study predicate $MV$ logic, then we are naturally led to expansions of $MV$ algebras allowing quantification.
\begin{definition} Let $\A$ be an $MV$ algebra. An existential quantifier on $\A$ is a function 
$\exists:A\to A$ that satisfies the following conditions for all $a, b\in A$:
\begin{enumerate}
\item $\exists 0=0.$
\item $a\leq \exists a.$
\item $\exists(a\odot \exists b)=\exists a\odot \exists b.$
\item $\exists(a\oplus \exists b)=\exists a\oplus \exists b.$
\item  $\exists (a\odot a)=\exists a\odot  \exists a.$
\item $\exists (a\oplus a)=\exists a\oplus \exists a.$
\end{enumerate}
\end{definition}
 Let $\A$ be an $MV$ algebra, with existential quantifier $\exists$. For $a\in A$, set $\forall a=\neg \exists \neg a$. 
Then $\forall$ is a unary operation on $\A$ called  a universal quantifier and it satisfies all properties
of the existential quantifier, except for (2) which takes the form: $\forall a\leq a$. This follows directly from the axioms.  
Now we define our algebras. Their similarity type depends on a fixed in advance semigroup. 
We write $X\subseteq_{\omega} Y$ to denote that $X$ is a finite subset 
of $Y$, and we write $\wp_{\omega}(Y)$ for $\{X: X\subseteq_{\omega}Y\}$.\footnote{There should be  no conflict with the notation $S_{\mathfrak{m}}V$ 
introduced earlier.}
\begin{definition} Let $\alpha$ be an infinite set. Let $G\subseteq {}^{\alpha}\alpha$ be a semigroup under the operation of composition of maps. Let 
$T\subseteq \wp(\alpha)$.
An $\alpha$ dimensional polyadic $MV$ algebra of type $(G,T)$, an $MV_{G,T}$ for short, 
is an algebra of the following
type
$$(A,\oplus,\odot,\neg, 0, 1, {\sf s}_{\tau}, {\sf c}_{(J)} )_{\tau\in G, J\in T}$$
where
$(A,\oplus,\odot, \neg, 0, 1)$ is an $MV$ algebra, ${\sf s}_{\tau}:\A\to \A$ is an endomorphism of $MV$ algebras,
${\sf c}_{(J)}$ is an existential quantifier, such that the following hold for all 
$p\in A$, $\sigma, \tau\in G$ and $J,J'\in T:$
\begin{enumarab}
\item ${\sf s}_{Id}p=p,$
\item ${\sf s}_{\sigma\circ \tau}p={\sf s}_{\sigma}{\sf s}_{\tau}p$ (so that $S:\tau\mapsto {\sf s}_{\tau}$ defines a homomorphism from $G$ to $End(\A)$; 
that is $(A, \oplus, \odot, \neg, 0, 1, G, S)$ is a transformation system),
\item ${\sf c}_{(J\cup J')}p={\sf c}_{(J)}{\sf c}_{(J')}p,$  
\item If $\sigma\upharpoonright \alpha\sim J=\tau\upharpoonright \alpha\sim J$, then
${\sf s}_{\sigma}{\sf c}_{(J)}p={\sf s}_{\tau}{\sf c}_{(J)}p,$ 
\item If $\sigma\upharpoonright \sigma^{-1}(J)$ is injective, then
${\sf c}_{(J)}{\sf s}_{\sigma}p={\sf s}_{\sigma}{\sf c}_{\sigma^{-1}(J)}p.$
\end{enumarab}
\end{definition}
\begin{theorem} Let $\A=(A,\oplus, \odot, \neg, 0, 1, {\sf s}_{\tau}, {\sf c}_{(J)})_{\tau\in G, J\in T}$ be an $MV_{G,T}$. 
For each $J\in T$ and $x\in A$, set ${\sf q}_{(J)}x=\neg {\sf c}_{(J)}\neg x.$
Then the following hold for each $p\in A$, $\sigma, \tau\in G$ and $J,J'\in T$.
\begin{enumarab}
\item ${\sf q}_{(J)}$ is a universal quantifier, 
\item ${\sf q}_{(J\cup J')}p={\sf q}_{(J)}{\sf c}_{(J')}p,$
\item ${\sf c}_{(J)}{\sf q}_{(J)}p={\sf q}_{(J)}p , \ \  {\sf q}_{(J)}{\sf c}_{(J)}p={\sf c}_{(J)}p,$
\item If $\sigma\upharpoonright \alpha\sim J=\tau\upharpoonright \alpha\sim J$, then
${\sf s}_{\sigma}{\sf q}_{(J)}p={\sf s}_{\tau}{\sf q}_{(J)}p,$
\item If $\sigma\upharpoonright \sigma^{-1}(J)$ is injective, then
${\sf q}_{(J)}{\sf s}_{\sigma}p={\sf s}_{\sigma}{\sf q}_{\sigma^{-1}(J)}p.$
\end{enumarab}
\end{theorem}
\begin{demo}{Proof} Routine
\end{demo}
Here we depart from \cite{HMT2} by defining polyadic algebras on sets rather than on ordinals. In this way we follow the tradition of Halmos.
We refer to $\alpha$ as the dimension of $\A$ and we write $\alpha=dim\A$.
Borrowing terminology from cylindric algebras, we refer to ${\sf c}_{(\{i\})}$ by ${\sf c}_i$ and ${\sf q}_{(\{i\})}$ by ${\sf q}_i$
\begin{example}
\begin{enumarab}
\item  Let $G$ be a semigroup of finite transformations on a set $V$ of variables, and let 
$T=\wp_{\omega}(V)$. Let $\mathfrak{L}_{\rho,G,T}$
be the corresponding logic. Assume that $V\sim \rho(p)$ is infinite for every $p\in P$. 
Let $Fm$ denote the set of formulas and $\Sigma\subseteq Fm$. 
Define on the set formulas the relation $\equiv_{\Sigma}$ by $\phi\equiv_{\Sigma}\psi$ if and only if 
$\Sigma\vdash \phi\longleftrightarrow \psi$ where $\longleftrightarrow$ is defined from $\to$ and $\odot$ the usual way.
Then the Tarski-Lindenbaum algeba
$\Fm/\equiv_{\Sigma}$, defined the obvious way, is a $G$ algebra.
The proof of this is tedious but straightforward. It says that the axioms and rules of inference of the logic $\mathfrak L_{G,T}$ 
and the polyadic axioms say the same thing under the proper interpretation. When $\rho(p)$ is finite for every $p\in P$ this is usual predicate
Lukasiewicz logic.

\item Let $G={}^VV$ and $T\in \{\wp(V), \wp_{\omega}(V)\}$. We need a modification to define the corresponding algebras. 
In case there are formulas $\phi$ such that $V_f(\phi)=V$ there is a problem to define the substitution operator
${\sf s}_{\tau}$. We did not encounter this difficulty in the previous example since 
we had an infinite supply of variables outside formulas, and those can be used to define the operation of simultaneous substitution.
We need extra variables in order to avoid collisions of free and bound variables. 
So let $V_1$ be a set of symbols which is disjoint from $P$ and $V$ 
and such that $|V_1|=|V|$. Let $V^*=V\cup V_1$, and let $\tau_0$ be a bijection between 
$V^*$ and $V_1$. Let $Fm$ be the set of formulas $\phi$ in the expanded language such that $V(\phi)\subseteq V$ for all $\phi\in Fm.$ 
Now $\Fm/\equiv_{\Sigma}$ is in  an $MV_{G,T}$ 
with substitutions defined as follows:
$${\sf s}_\tau(\phi/\Sigma)=({\sf S}_f(\tau){\sf S}_f(\tau_0^{-1}){\sf S}({\tau_0}))\phi/\Sigma.$$
(The other operations are defined the usual way).
\item Let $L$ be a fixed linearly ordered complete $MV$ algebra. Let $X$ and $I$ be any sets. 
For $x,y\in {}^IX$ and $J\subseteq I$, write $x\equiv_J y$ if $x(l)=y(l)$ for all $l\notin J$. For $p: {}^IX\to L$, we let
$${\sf c}_{(J)}p(x)=\bigvee\{p(y): x\equiv_J y\},$$
and $${\sf s}_{\tau}p(x)=p(x\circ \tau).$$
The typical example is when $L=[0,1]$. Such special algebras will be denoted by $\F(^IX, [0,1]).$ If $V\subseteq X^I$, 
we let $\F(V,[0,1])$ denote the algebra with universe all function from $V$ to $[0,1]$ with operations as above.
Not all such algebras are $MV$ polyadic algebras, but for some choices of $V,$ 
the resulting algebra is an $MV$ polyadic algebra, see theorem \ref{weak}. 

\item Let $G={}^VV$ and $T\subseteq \wp(V).$ Let $\M=(M,p_M)_{p\in P}$ be a model. 
For a formula $\phi$,  let $\phi^{\M}=\{s: {}^VM\to [0,1]: ||\phi||_{\M,s}=1\}$. Let $Fm$ denote the 
set of formulas. 
Then $\{\phi^{\M}: \phi\in Fm\}$ is the universe of a polyadic $MV$ algebra, which we denote by $\C^{\M}$ to emphasize the role played by $\M$ 
and call it the set algebra based on $\M$.
The operations are read off from the semantics of connectives. For example
$$\phi^{\M}\oplus \psi^{\M}=(\phi\oplus \psi)^{\M},$$
and $${\sf c}_{(J)}\phi^{\M}=(\exists J\phi)^{\M}.$$  It is easy to observe that such algebras are special cases from the algebras in the preceding 
item, when $L=[0,1]$.
\end{enumarab}
\end{example}
\begin{definition}\label{rep} An $MV$ algebra is representable, if for all $a\neq 0$, there exist a set $X$, $V\subseteq {}^IX$ and a homomorphism 
from $\A$ to $\F[V,[0,1])$
such that  $f(a)\neq 0$. 
\end{definition}
$MV$ algebras of the form $\F(V,[0,1])$ are called set algebras. 
We define $V$ that gives rise to such ``concrete set algebras":
\begin{theorem}\label{weak} For two given sets $I$ and $X$ and $p:I\to X$, let $^IX^{(p)}$ be the following set
$\{s\in {}^IX: |\{i\in I: s_i\neq p_i\}|<\omega\}.$ Such a set is called a weak space. For a union of weak spaces $W$, 
$F(W, [0,1])$ is the universe of 
an $MV$ polyadic algebra with operations defined as above.
\end{theorem}
\begin{demo}{Proof} Direct.
\end{demo}

For the case $G={}^{\alpha}\alpha$ and $T=\wp_{\omega}(\alpha)$, 
we have a natural correspondence between equations in algebras and equivalences in ${\mathfrak L}_{G,T}$.
We formulate this correspondence for $\alpha=\omega$, and special languages where 
predicates are countable, and atomic formulas are of of the form $p(v_0, v_1\ldots v_i\ldots)_{i<\omega}$, 
that is, we have countably many predicates each of arity $\omega$ and in atomic formulas 
variables occur only in their natural order. This is not really any different 
from usual languages when arity of atomic formulas is $\omega$, since other atomic formulas can be recovered by applying 
substitutions to the restricted ones. A similar correspondence can be obtained when $T=\wp(\alpha)$.
The same can be done for strongly rich semigroups on $\omega$ (to be defined shortly), but
for $G$ consisting of finite transformations it is more involved, cf. \cite{HMT2} theorems 4.3.58, 4.3.59.
\begin{definition} Let $\alpha$ be an infinite set. $L_{\alpha}$ is the standard language of the class $MV_{G,T}$; 
it has $\omega$ variables $x_0,x_1\ldots$ and operation symbols 
$\oplus, \odot, \neg, 0, 1, {\sf c}_k, {\sf s}_{\tau},$ with $k\in \alpha$ and $\tau\in {}^\alpha\alpha$. Let $\Lambda=(V, P, G)$ be a language with $G={}^VV$, 
$|P|=\omega$; 
we assume that
$P=\{p_0, p_1\ldots p_i\ldots\ :i\in \omega\}.$
With each term $\sigma$ of $L_{\alpha},$ 
we associate a formula $\eta\sigma$ of $\Lambda$, with $i, k\in \omega$, as follows: 
$$\eta x_i=P_i(v_0,\ldots v_i\ldots )_{i<\omega},$$
$$\eta(\sigma\oplus \tau)=(\eta\sigma\oplus \eta\tau),$$
$$\eta(\sigma\odot \tau)=(\eta\sigma\odot \eta\tau),$$
$$\eta(0)=\bot,$$
$$\eta(1)=\top,$$
$$\eta(\neg \sigma)=(\neg \sigma),$$
$$\eta ({\sf c}_k\sigma) =(\exists v_k\eta \sigma),$$
$$\eta ({\sf s}_{\tau}\sigma)={\sf S}(\tau)(\eta\sigma).$$
\end{definition}
The following is a typical translation that abounds in algebraic logic, cf. \cite{HMT2} theorem 4.3.57.
\begin{theorem}\label{model}
\begin{enumarab}
\item Let $\M$ be a a structure  
and $\A=\C^{\M}$.  Then for any term $\sigma$
$$(\eta(\sigma))^{\M}=(\sigma^{\A}P).$$
Here $P$ is looked at as an assignment $P:\omega\to \A$, such that  
$$P_i=p_i^{\M}=\{s:V\to [0,1]: ||p_i||_{s,\M}=1\}.$$
\item For any $\sigma, \tau$ of $L_{\alpha}$ the following conditions are equivalent
\begin{enumroman}
\item $\models \eta\sigma\leftrightarrow \eta\tau$
\item $MV_{G,T}\models \sigma=\tau.$
\end{enumroman}
\end{enumarab}
\end{theorem}
\begin{demo}{Sketch of Proof} We prove only $(ii)\to (i)$. Let $\A$ be a set algebra with universe $A=F(^IX, [0,1])$. Let $a\in {}^{\omega}A.$
Let $\D$ be the subalgebra generated by the range of $a$.  Define a model $\M$, by stipulating that for $s\in {}^IX$, 
$|||R_i[s]||=a_i(s).$ Then $\C^{\M}= \D$. 
Now $(\eta\sigma)^{\M}=(\eta\tau)^{\M}$, so $\sigma a=\tau a$, and we are done.

\end{demo}
We shall deal with the case when $G$ is a countable special proper subsemigroup of $^{\alpha}\alpha$.
Here $\alpha$ is a countable set (sometimes it will be an ordinal) and algebras considered are also countable.
We need some preparations to define such semigroups.

\begin{athm}{Notation} For  a set $X$, recall that $|X|$ stands for the cardinality of $X$.
$Id_X$, or simply $Id$ when $X$ is clear from context,  denotes the identity function
on $X$.
For functions $f$ and $g$ and a set $H$, $f[H|g]$ is the function
that agrees with $g$ on $H$, and is otherwise equal to $f$. $Rgf$ denotes 
the range of $f$. For a transformation $\tau$ on $\omega$, 
the support of $\tau$, or $sup(\tau)$ for short, is the set:
$$sup(\tau)=\{i\in \omega: \tau(i)\neq i\}.$$
Let $i,j\in \omega$, then $\tau[i|j]$ is the transformation on $\omega$ defined as follows:
$$\tau[i|j](x)=\tau(x)\text { if } x\neq i \text { and }\tau[i|j](i)=j.$$ 
For a function $f$, $f^n$ denotes the composition $f\circ f\ldots \circ f$
$n$ times.
\end{athm}

\begin{definition}  Let $\alpha$ be a countable set. Let $T\subseteq \langle {}^{\alpha}\alpha, \circ \rangle$ be a semigroup.
We say that $T$ is {\it rich } if $T$ satisfies the following conditions:
\begin{enumerate}
\item $(\forall i,j\in \alpha)(\forall \tau\in T) \tau[i|j]\in T.$
\item There exists $\sigma,\pi\in T$ such that
$(\pi\circ \sigma=Id,\  Rg\sigma\neq \alpha).$
\item $ (\forall \tau\in T)(\sigma\circ \tau\circ \pi)[(\alpha\sim Rg\sigma)|Id]\in T.$
\item Let $T\subseteq \langle {}^{\alpha}\alpha, \circ\rangle$ be a rich semigroup. 
Let $\sigma$ and $\pi$ be as in the previous item. 
If $\sigma$ and $\pi$ satisfy $(i)$, $(ii)$ below:
$$ (i)\ \ \ \  (\forall n\in \alpha) |supp(\sigma^n\circ \pi^n)|<\alpha.$$  
$$(ii) \ \ \ \  (\forall n\in \alpha)[supp(\sigma^n\circ \pi^n)\subseteq 
\alpha\smallsetminus Rg(\sigma^n)];$$
then we say that $T$ is  {\it a strongly rich} semigroup. 
\end{enumerate}
\end{definition}

\begin{example}The semigroup 
$(^{\omega}\omega, \circ)$ is rich (but not countable) and so is its  subsemigroup generated by
$\{[i|j], [i,j], suc, pred\}$. Here $suc$ abbreviates the successor function on $\omega$ and $pred$ is the function defined by 
$pred(0)=0$ and for other $n\in \omega$, $pred(n)=n-1$. In fact, both semigroups are strongly rich;
in the second case $suc$ plays the role of $\sigma$ while $pred$ plays the role of 
$\pi$.
\end{example}
Next, we collect some properties of $MV_{G,T}$ algebras that are more handy to use in our subsequent work.

\begin{theorem}\label{d} Let $\alpha$ be an infinite set, $G\subseteq {}^{\alpha}\alpha$ a semigroup, $T=\wp_{\omega}(\alpha)$, and $\A\in MV_{G,T}$. 
Then $\A$ satisfies the following for $\tau,\sigma\in G$
and all $i,j,k\in \alpha$.
\begin{enumerate}

\item $x\leq {\sf c}_ix={\sf c}_i{\sf c}_ix,\ {\sf c}_i(x\oplus c_iy)={\sf c}_ix\oplus {\sf c}_iy,\ {\sf c}_i(-{\sf c}_ix)=-{\sf c}_ix,\ 
{\sf c}_i{\sf c}_jx={\sf c}_j{\sf c}_ix$.

\item ${\sf s}_{\tau}$ is an $MV$ algebra  endomorphism.

\item  ${\sf s}_{\tau}{\sf s}_{\sigma}x={\sf s}_{\tau\circ \sigma}x$
and ${\sf s}_{Id}x=x$.

\item ${\sf s}_{\tau}{\sf c}_ix={\sf s}_{\tau[i|j]}{\sf c}_ix$.

Recall that $\tau[i|j]$ is the transformation that agrees with $\tau$ on 
$\alpha\smallsetminus\{i\}$ and $\tau[i|j](i)=j$.

\item ${\sf s}_{\tau}{\sf c}_ix={\sf c}_j{\sf s}_{\tau}x$ if $\tau^{-1}(j)=\{i\}$, 
${\sf s}_{\tau}{\sf q}_ix={\sf q}_j{\sf s}_{\tau}x$ 
if $\tau^{-1}(j)=\{i\}$.  

\item  ${\sf c}_i{\sf s}_{[i|j]}x={\sf s}_{[i|j]}x$,\ \ ${\sf q}_i{\sf s}_{[i|j]}x={\sf s}_{[i|j]}x.$

\item  ${\sf s}_{[i|j]}{\sf c}_ix={\sf c}_ix$, \ \ ${\sf s}_{[i|j]}{\sf q}_ix={\sf q}_ix$.

\item ${\sf s}_{[i|j]}{\sf c}_kx={\sf c}_k{\sf s}_{[i|j]}x$,\ \ ${\sf s}_{[i|j]}{\sf q}_kx={\sf q}_k{\sf s}_{[i|j]}x$
whenever $k\notin \{i,j\}$.

\item  ${\sf c}_i{\sf s}_{[j|i]}x={\sf c}_j{\sf s}_{[i|j]}x$,\ \  ${\sf q}_i{\sf s}_{[j|i]}x={\sf q}_j{\sf s}_{[i|j]}x$.
\end{enumerate}
\end{theorem}
\begin{demo}{Proof} The proof is tedious but fairly straighforward.
\end{demo}
Following cylindric algebra terminology, we will be often writing ${\sf s}_j^i$ for ${\sf s}_{[i|j]}$.
We now define an important concept in algebraic logic. This concept occurs under the rubric  of neat reducts
in cylindric algebras \cite{HMT1}, \cite{AUneat} and compressions in polyadic algebras \cite{DM}.
\begin{definition} 
\begin{enumarab}
\item Let $\alpha\subseteq \beta$ be infinite sets, and let $G_{\beta}$ be a semigroup of transformations on $\beta$ and 
$G_{\alpha}$ be a semigroup of transformations on $\alpha$.
For $\tau$ in $G_{\alpha}$ we write $\bar{\tau}$ for $\tau\cup Id_{\beta\sim \alpha}.$ 
\item  Let $ \alpha\subseteq \beta$ be infinite sets. Let $G_{\beta}$ be a semigroup of transformations on $\beta$, 
and let $G_{\alpha}$ be a semigroup of transformations on $\alpha$ such that for all $\tau\in G_{\alpha}$, one has $\bar{\tau}\in G_{\beta}$. 
Let $T\subseteq \wp(\beta)$.
Let $\A=(A, \oplus, \odot, \neg, 0,1, {\sf c}_{(J)}, {\sf s}_{\tau})_{J\in T, \tau\in G_{\beta}}$ be an algebra in $MV_{G_{\beta},T}$.
Then $\Rd_{\alpha}\A$ is the $G_{\alpha}$ algebra obtained by discarding operations indexed by elements in  $\beta\sim \alpha$. That is 
$\Rd_{\alpha}\A=(A, \oplus, \odot, \neg,  0,1, {\sf c}_{(J)}, {\sf s}_{\bar{\tau}})_{J\in T\cap \wp(\alpha), \tau\in G_{\alpha}}$.
\item Let $T=\wp_{\omega}(\beta)$. For $\A\in MV_{G_{\beta},T}$  
and $x\in A$, then $\Delta x,$ the dimension set of $x$, 
is defined by $\Delta x=\{i\in \beta: {\sf c}_ix\neq x\}.$
Let $B=\{x\in A: \Delta x\subseteq \alpha\}$. If $B$ happens to be 
a subuniverse of $\Rd_{\alpha}\A$, 
then $\B$ is a subreduct of $\A$, it is  called the $\alpha$ neat reduct of $\A$ and is denoted by $\Nr_{\alpha}\A$.
\item For $T=\wp(\beta)$ and $\A\in MV_{G_{\beta},T}$,  let $B=\{x\in A: {\sf c}_{(\beta\sim \alpha)}x=x\}$. 
If $B$ happens to be a subuniverse of $\Rd_{\alpha}\A$, then it is called the $\alpha$ neat reduct of $\A$ and is denoted by $\Nr_{\alpha}\A$.
\item Let $K$ be a class of algebras having dimension $\beta$ and $\alpha\subseteq \beta$. 
Assume that the operator $\Nr_{\alpha}$ makes sense, that is, it is applicable to $K$.
Then $\Nr_{\alpha}K=\{\Nr_{\alpha}\A: \A\in K\}.$ 

\end{enumarab}
\end{definition}
For an algebra $\A$, and $X\subseteq \A$, $\Sg^{\A}X$ or simply $\Sg X$, when $\A$ is clear from context,
denotes the subalgebra of $\A$ generated by $X$
\begin{theorem}\label{rich} Let $\alpha\subseteq \beta$ be countably infinite sets. If $G$ is a strongly rich semigroup on $\alpha$ 
and $\A\in MV_{G,T}$, where $T=\wp_{\omega}(\alpha)$, then there exists a strongly rich semigroup $H$ on $\beta$ 
and $\B\in MV_{H, \bar{T}}$, where $\bar{T}=\wp_{\omega}(\beta)$ such that
$\A\subseteq \Nr_{\alpha}\B$ and for all $X\subseteq A,$ one has  $\Sg^{\A}X=\Nr_{\alpha}\Sg^{\B}X$. 
\end{theorem}

\begin{demo}{Proof}  cf. \cite{AUamal}. We assume that $\alpha$ is an ordinal; in fact without loss of generality we can assume that it 
is the least infinite ordinal $\omega.$
We also assume a particular strongly rich semigroup, namely that generated by finite transformations together with 
$suc$, $pred$. The general case is the same \cite{AUamal} Remark 2.8 p.327.
We follow \cite {AUamal} p.323-336, referring to op.cit for detailed arguments. 
Let $n\leq \omega$. Then $\alpha_n=\omega+n$ 
and $M_n=\alpha_n\sim \omega$.
Note that when $n\in \omega$, then $M_n=\{\omega,\ldots,\omega+n-1\}$.
Let $\tau\in G$, then $\tau_n=\tau\cup Id_{M_n}$. $T_n$ denotes the 
subsemigroup of $\langle {}^{\alpha_n}\alpha_n,\circ \rangle$ generated by
$\{\tau_n:\tau\in G\} \cup \cup_{i,j\in \alpha_n}\{[i|j],[i,j]\}$.
For $n\in \omega$, let $\rho_n:\alpha_n\to \omega$ 
be the bijection defined by 
$\rho_n\upharpoonright \omega=suc^n$ and $\rho_n(\omega+i)=i$ for all $i<n$.
Let $n\in \omega$. For $v\in T_n,$ let $v'=\rho_n\circ v\circ \rho_n^{-1}$.
Then $v'\in G$.
For $\tau\in T_{\omega}$, let 
$D_{\tau}=\{m\in M_{\omega}:\tau^{-1}(m)=\{m\}=\{\tau(m)\}\}$.
Then $|M_{\omega}\sim  D_{\tau}|<\omega.$ 
Let $\A$ is an  arbitrary countable $G$ algebra.  
Let $\A_n$ be the algebra defined as follows:
$\A_n=\langle A,\oplus, \odot, \neg, 0, 1, {\sf c}_i^{\A_n},{\sf s}_v^{\A_n}\rangle_{i\in \alpha_n,v\in T_n}$
where for each $i\in \alpha_n$ and $v\in T_n$, 
${\sf c}_i^{\A_n}:= {\sf c}_{\rho_n(i)}^{\A} \text { and }{\sf s}_v^{\A_n}:= {\sf s}_{v'}^{\A}.$
Let $\Rd_{\omega}\A_n$ be the following reduct of 
$\A_n$ obtained by restricting the type of $\A_n$ to the first 
$\omega$ dimensions:
$\Rd_{\omega}\A_n=\langle A_n,\oplus,\odot, \neg, 0, 1,  {\sf c}_i^{\A_n}, {\sf s}_{\tau_n}^{\A_n}\rangle_{i\in \omega,\tau\in G}.$
For $x\in A$, let $e_n(x)={\sf s}_{suc^n}^{\A}(x)$. Then $e_n:A\to A_n$ and  $e_n$ is an isomorphism 
from $\A$ into $\Rd_{\omega}\A_n$
such that $e_n(\Sg^{\A}Y)=\Nr_{\omega}(\Sg^{\A_n}e_n(Y))$ for all
$Y\subseteq A$, cf. \cite{AUamal} claim 2.7.\footnote{ In \cite{AUamal} the algebras have a Boolean reduct. But such maps 
preserve the $MV$ operations, in our present context,  because the substitutions are endomorphisms of the algebra in question.}
For the sake of brevity, let $\alpha=\alpha_{\omega}=\omega+\omega$.
Let $T_{\omega}$ is the semigroup generated by the set
$\{\tau_{\omega}: \tau\in G\}\cup_{i,j\in \alpha}\{[i|j],[i,j]\}.$
For $\sigma\in T_{\omega}$, and $n\in \omega$, 
let $[\sigma]_n=\sigma\upharpoonright \omega+n$.
For each $n\in \omega,$ let 
$\A_n^+=\langle A,\oplus,\odot, \neg, 0, 1, {\sf c}_i^{\A_n^+}, {\sf s}_{\sigma}^{\A_n^+}\rangle_{i\in \alpha, \sigma\in T_{\omega}}$
be an expansion of $\A_n$ 
such that there $MV$ reducts coincide  
and for each 
$\sigma\in T_{\omega}$ and $i\in \alpha,$ 
${\sf s}_{\sigma}^{\A_n^+}:={\sf s}_{[\sigma]_n}^{\A_n} 
\text { iff } [\sigma]_n\in T_n,$ 
and 
${\sf c}_i^{\A_n^+}:={\sf c}_i^{\A_n}\text { iff }i<\omega+n.$
Let $F$ be any non-principal ultrafilter on $\omega$. Now 
forming the ultraproduct of the $\A_n^+$'s relative to $F$, let
$\A^+=\prod_{n\in \omega}\A_n^+/F.$
For $x\in A$, let 
$e(x)=\langle e_n(x):n\in \omega\rangle/F.$
Let $\Rd_{\omega}A^+=\langle A^+, \oplus, \odot, \neg, 0, 1, {\sf c}_i^{\A^+}, {\sf s}_{\tau_{\omega}}^{\A^+}
\rangle_{i<\omega,\tau\in T}.$
Then $e$ is an isomorphism from $\A$ into $\Rd_{\omega}\A^+$
such that  $e(\Sg^{\A}Y)=\Nr_{\omega}\Sg^{\A^+}e(Y)$ for all $Y\subseteq A.$
We have shown that $\A$ neatly embeds in algebras in finite extra dimensions and in $\omega$ extra dimension. An iteration of this embedding yields the required result. 
\end{demo}
 
\begin{theorem}\label{net}
Let $\alpha\subseteq \beta$ be infinite sets. Then the following hold:
\begin{enumarab}
\item Let $G_{I}$ be the semigroup of finite transformations on $I$.
Let $\A\in MV_{G_{\alpha},T}$ where $T=\wp_{\omega}(\alpha)$ and $\alpha\sim \Delta x$ is infinite for every $x\in A$. 
Then there exists $\B\in MV_{G_{\beta},\bar{T}}$, where $\bar{T}=\wp_{\omega}(\beta)$ 
such that $\A\subseteq \Nr_{\alpha}\B$ and for all $X\subseteq 
A$, one has $\Sg^{\A}X=\Nr_{\alpha}\Sg^{\B}X.$

\item  Let $G_I$ be the semigroup of all transformations on $I$. Let $T=\wp(\alpha)$.
Let $\A\in MV_{G_{\alpha},T}.$
Then for all $\beta\supseteq \alpha$, there exists $\B\in MV_{G_{\beta}, \bar{T}}$, where $\bar{T}=\wp(\beta)$, 
such that $\A\subseteq \Nr_{\alpha}\B$ and for all $X\subseteq 
A$ one has $\Sg^{\A}X=\Nr_{\alpha}\Sg^{\B}X.$ A completely analagous result holds  by replacing $T$ and $\bar{T}$ by $\wp_{\omega}(\alpha)$ and 
$\wp_{\omega}(\beta)$, respectively.
\end{enumarab}
\end{theorem}
\begin{demo}{Proof}
\begin{enumarab}
\item  Let $\alpha\subseteq \beta$. We assume that $\alpha$ and $\beta$ are ordinals with $\alpha<\beta$.
The proof is an adaptation of the proof of Theorem  2.6.49 (i)
in \cite{HMT1}. First we show that there exists $\B\in MV_{G_{\alpha+1},\bar{T}}$, where $\bar{T}=\wp_{\omega}(\alpha+1)$ 
such that $\A$ embeds into $\Nr_{\alpha}\B$. 
Let $$R = Id\upharpoonright (\alpha\times A)
 \cup \{  ((k,x), (\lambda, y)) : k, \lambda <
\alpha, x, y \in A, \lambda \notin \Delta x, y = {\mathsf s}_{[k|\lambda]} x \}.$$
It is easy to see  that $R$ is an equivalence relation on $\alpha
\times A$.
Define the following operations on $(\alpha\times A)/R$ with $\mu, i, k\in \alpha$ and $x,y\in A$:
(Abusing notation we denote the new operations like the old ones. No confusion should ensue).
\begin{equation*}\label{l5}
\begin{split}
(\mu, x)/R \oplus (\mu, y)/R = (\mu, x \oplus y)/R, 
\end{split}
\end{equation*}
\begin{equation*}\label{l6}
\begin{split}
(\mu, x)/R\odot  (\mu, y)/R = (\mu, x\odot  y)/R, 
\end{split}
\end{equation*}
 \begin{equation*}\label{l7}
\begin{split}
\neg (\mu, x)/R=(\neg \mu, x)/R,
\end{split}
\end{equation*}
\begin{equation*}\label{l10}
\begin{split}
{\mathsf c}_i ((\mu, x)/R)  = (\mu, {\mathsf c}_i x )/R, \quad
\mu \in \alpha \smallsetminus
\{i\},
\end{split}
\end{equation*}
\begin{equation*}\label{l11}
\begin{split}
{\mathsf s}_{[j|i]} ((\mu, x)/R)  = (\mu, {\mathsf s}_{[j|i]} x )/R, \quad \mu \in \alpha
\smallsetminus \{i, j\}.
\end{split}
\end{equation*}
The constants $0$ and $1$ are defined in the obvious way.
It can be checked that these operations are well defined.
Let $$\C=((\alpha\times A)/R, \oplus, \odot, \neg, 0, 1,  {\sf c_i}, {\sf s}_{i|j]})_{i,j\in \alpha},$$
be the quotient algebra and
let $$h=\{(x, (\mu,x)/R): x\in A, \mu\in \alpha\sim \Delta x\}.$$
Then $h$ is well defined, indeed $h$ is an isomorphism from $\A$ into $\C$. 
Now to show that $\A$ neatly embeds into $\alpha+1$ extra dimensions we define the operations ${\sf c}_{\alpha}, {\sf s}_{[i|\alpha]}$
and ${\sf s}_{[\alpha|i]}$ on $\C$ as follows:
$${\mathsf c}_\alpha = \{ ((\mu, x)/R, (\mu, {\mathsf c}_\mu x)/R) :
\mu \in \alpha, x \in B \},$$
$${\mathsf s}_{[i|\alpha]} = \{ ((\mu, x)/R, (\mu, {\mathsf s}_{[i|\mu]}
x)/R) : \mu \in \alpha \smallsetminus \{i\}, x \in B \},$$
$${\mathsf s}_{[\alpha|i]} = \{ ((\mu, x)/R, (\mu, {\mathsf s}_{[\mu|i]}
x)/R) : \mu \in \alpha \smallsetminus \{i\}, x \in B \}.$$
Let $$\B=((\alpha\times A)/R, \oplus,\odot, \neg, 0, 1, {\sf c}_i, {\sf s}_{[i|j]})_{i,j\leq \alpha}.$$
Then $\B\in MV_{G_{\alpha+1},\bar{T}}$ where $\bar{T}=\wp_{\omega}(\alpha+1) \text{ and }h(\A)\subseteq \Nr_{\alpha}\B.$
It is not hard to check that the defined operations are as desired. We have our result when $G$ consists only of replacements.
But since $\alpha\sim \Delta x$ is infinite one can show that substitutions corresponding to all finite transformations are term definable as follows.
For a given finite transformation $\tau\in {}^{\alpha}\alpha,$ we write $[u_0|v_0, u_1|v_1,\ldots,
u_{k-1}|v_{k-1}]$ if $sup(\tau)=\{u_0,\ldots ,u_{k-1}\}$, $u_0<u_1
\ldots <u_{k-1}$ and $\tau(u_i)=v_i$ for $i<k$.
Let $\A\in MV_{G,T}$ be such that $\alpha\sim \Delta x$ is
infinite for every $x\in A$. If $\tau=[u_0|v_0, u_1|v_1,\ldots,
u_{k-1}|v_{k-1}]$ is a finite transformation, if $x\in A$ and if
$\pi_0,\ldots ,\pi_{k-1}$ are in this order the first $k$ ordinals
in $\alpha\sim (\Delta x\cup Rg(u)\cup Rg(v))$, then
$${\mathsf s}_{\tau}x={\mathsf s}_{v_0}^{\pi_0}\ldots
{\mathsf s}_{v_{k-1}}^{\pi_{k-1}}{\mathsf s}_{\pi_0}^{u_0}\ldots
{\mathsf s}_{\pi_{k-1}}^{u_{k-1}}x.$$
The ${\sf s}_{\tau}$'s so defined satisfy the polyadic axioms, cf. \cite{HMT1} Theorem 1.11.11.
Then one proceeds by a simple induction to show that for all $n\in \omega$ there exists $\B\in MV_{G_{\alpha+n},\bar{T}}$ where 
$\bar{T}=\wp_{\omega}(\alpha+n)$ 
such that $\A\subseteq \Nr_{\alpha}\B.$ For the transfinite, one uses ultraproducts, cf. \cite{HMT1} theorem 2.6.34. 
For the second part, let $\A\subseteq \Nr_{\alpha}\B$ and $A$ generates $\B,$ then $\B$ consists of all elements ${\sf s}_{\sigma}^{\B}x$ such that 
$x\in A$ and $\sigma$ is a finite transformation on $\beta$ such that
$\sigma\upharpoonright \alpha$ is one to one, cf. \cite{HMT1} lemma 2.6.66. 
Now suppose $x\in \Nr_{\alpha}\Sg^{\B}X$ and $\Delta x\subseteq
\alpha$. There exist $y\in \Sg^{\A}X$ and a finite transformation $\sigma$
of $\beta$ such that $\sigma\upharpoonright \alpha$ is one to one
and $x={\sf s}_{\sigma}^{\B}y.$  
Let $\tau$ be a finite
transformation of $\beta$ such that $\tau\upharpoonright  \alpha=Id
\text { and } (\tau\circ \sigma) \alpha\subseteq \alpha.$ Then
$x={\sf s}_{\tau}^{\B}x={\sf s}_{\tau}^{\B}{\sf s}_{\sigma}y=
{\sf s}_{\tau\circ \sigma}^{\B}y={\sf s}_{\tau\circ
\sigma\upharpoonright \alpha}^{\A}y.$

\item Let $(\A, \alpha,S)$ be a transformation system. 
That is $\A$ is an $MV$ algebra and $S:({}^\alpha\alpha,\circ)\to End(\A)$ is a homomorphism of semigroups. 
For any set $X$, let $F(^{\alpha}X,\A)$ 
be the set of all functions from $^{\alpha}X$ to $\A$ endowed with $MV$ operations defined pointwise and for 
$\tau\in {}^\alpha\alpha$ and $f\in F(^{\alpha}X, \A)$, ${\sf s}_{\tau}f(x)=f(x\circ \tau)$. (These are cylindrification free reducts of set algebras).
This turns $F(^{\alpha}X,\A)$ to a transformation system as well. 
The map $H:\A\to F(^{\alpha}\alpha, \A)$ defined by $H(p)(x)={\sf s}_xp$ is
easily checked to be an isomorphism. Assume that $\beta\supseteq \alpha$. Then $K:F(^{\alpha}\alpha, \A)\to F(^{\beta}\alpha, \A)$ 
defined by $K(f)x=f(x\upharpoonright \alpha)$ is an isomorphism. These facts are straighforward to establish, cf. theorem 3.1, 3.2 
in \cite{DM}. 
$F(^{\beta}\alpha, \A)$ is called a minimal dilation of $F(^{\alpha}\alpha, \A)$ of dimension $\beta$. 
Elements of the big algebra, or the $\beta$ - dilation, are of form ${\sf s}_{\sigma}p$,
$p\in F(^{\beta}\alpha, \A)$ where $\sigma$ is one to one on $\alpha$, cf. \cite{DM} theorem 4.3-4.4.
We say that $J\subseteq I$ supports an element $p\in A$ if whenever $\sigma_1$ and  $\sigma_2$ are 
transformations that agree on $J,$ then  ${\sf s}_{\sigma_1}p={\sf s}_{\sigma_2}p$.
$\Nr_JA$, consisting of the elements that $J$ supports, is called a compression of $\A$; 
with the operations defined the obvious way. (Note that this is the cylindrifier free definition of neat reducts).
If $\A$ is a $\B$ valued $I$ transformaton system wih domain $X$, 
then the $J$ compression of $\A$ is isomorphic to a $\B$ valued $J$ transformation system
via $H: \Nr_J\A\to F(^JX, \A)$ by setting for $f\in\Nr_J\A$ and $x\in {}^JX$, $H(f)x=f(y)$ where $y\in X^I$ and $y\upharpoonright J=x$, 
cf. \cite{DM} theorem 3.10.
Now let $\alpha\subseteq \beta.$ If $|\alpha|=|\beta|$ then the the required algebra is defined as follows. 
Let $\mu$ be a bijection from $\beta$ onto $\alpha$. For $\tau\in {}^{\beta}\beta,$ let ${\sf s}_{\tau}={\sf s}_{\mu\tau\mu^{-1}}$ 
and for each $i\in \beta,$ let 
${\sf c}_i={\sf c}_{\mu(i)}$. Then this defines $\B\in MV_{G,T}$, where $T=\wp(\beta)$, in which $\A$ neatly embeds via 
${\sf s}_{\mu\upharpoonright\alpha},$ cf.
\cite{DM} p.138.  Now assume that $|\alpha|<|\beta|$ .
Let $\A$ be a  given polyadic algebra of dimension $\alpha$; discard its cylindrifiers and then take its minimal dilation $\B$ of dimension $\beta$, 
which exists by the above. 
We need to define cylindrifiers on the big algebra, so that they agree with their values in $\A$ and to have $\A\cong \Nr_{\alpha}\B$. 
For $T=\wp(\alpha)$, we set, see \cite{DM} p. 165,
$${\sf c}_{(\Gamma)}{\sf s}_{\sigma}^{\B}p={\sf s}_{\rho^{-1}}^{\B} 
{\sf c}_{(\rho(\Gamma)\cap \sigma \alpha)}{\sf s}_{(\rho\sigma\upharpoonright \alpha)}^{\A}p,$$
and for $T=\wp_{\omega}(\alpha)$, we set (*):
$${\sf c}_k{\sf s}_{\sigma}^{\B}p={\sf s}_{\rho^{-1}}^{\B} {\sf c}_{(\rho\{k\}\cap \sigma \alpha)}{\sf s}_{(\rho\sigma\upharpoonright \alpha)}^{\A}p.$$
In the last two equations  $\rho$ is a any permutation on $\beta$. It can be checked that this definition is sound; it is independent of the choice of $\rho$. 
Furthermore, it defines the required algebra $\B$. Let us check this.
Since our definition is different from that in \cite{DM}, by restricting cylindrifiers to be only finite, and the algebras considered are expansions of 
$MV$ algebras rather than Boolean algebras,  we need to check the polyadic axioms which is tedious but basically routine.
We check only the axiom $${\sf c}_k(q_1\odot {\sf c}_kq_2)={\sf c}_kq_1\odot {\sf c}_kq_2.$$
All the same, we follow closely \cite{DM} p. 166. 
Assume that $q_1={\sf s}_{\sigma}^{\B}p_1$ and $q_2={\sf s}_{\sigma}^{\B}p_2$. 
Let $\rho$ be a permutation of $I$ such that $\rho(\sigma_1I\cup \sigma_2I)\subseteq I$ and let 
$$p={\sf s}_{\rho}^{\B}[q_1\odot {\sf c}_kq_2].$$
Then $$p={\sf s}_{\rho}^{\B}q_1\odot {\sf s}_{\rho}^{\B}{\sf c}_kq_2
={\sf s}_{\rho}^{\B}{\sf s}_{\sigma_1}^{\B}p_1\odot {\sf s}_{\rho}^{\B}{\sf c}_k {\sf s}_{\sigma_2}^{\B}p_2.$$
Now we calculate ${\sf c}_k{\sf s}_{\sigma_2}^{\B}p_2.$
We have by (*)
$${\sf c}_k{\sf s}_{\sigma_2}^{\B}p_2= {\sf s}^{\B}_{\sigma_2^{-1}}{\sf c}_{\rho(\{k\}\cap \sigma_2I)} {\sf s}^{\A}_{(\rho\sigma_2\upharpoonright I)}p_2.$$
Hence $$p={\sf s}_{\rho}^{\B}{\sf s}_{\sigma_1}^{\B}p_1\odot {\sf s}_{\rho}^{\B}{\sf s}^{\B}_{\sigma^{-1}}{\sf c}_{\rho(\{k\}\cap \sigma_2I)} 
{\sf s}^{\A}_{(\rho\sigma_2\upharpoonright I)}p_2.$$
\begin{equation*}
\begin{split}
&={\sf s}^{\A}_{\rho\sigma_1\upharpoonright I}p_1\odot {\sf s}_{\rho}^{\B}{\sf s}^{\A}_{\sigma^{-1}}{\sf c}_{\rho(\{k\}\cap \sigma_2I)} 
{\sf s}^{\A}_{(\rho\sigma_2\upharpoonright I)}p_2,\\
&={\sf s}^{\A}_{\rho\sigma_1\upharpoonright I}p_1\odot {\sf s}_{\rho\sigma^{-1}}^{\A}
{\sf c}_{\rho(\{k\}\cap \sigma_2I)} {\sf s}^{\A}_{(\rho\sigma_2\upharpoonright I)}p_2,\\
&={\sf s}^{\A}_{\rho\sigma_1\upharpoonright I}p_1\odot {\sf c}_{\rho(\{k\}\cap \sigma_2I)} {\sf s}^{\A}_{(\rho\sigma_2\upharpoonright I)}p_2.\\
\end{split}
\end{equation*} 
Now $${\sf c}_k{\sf s}_{\rho^{-1}}^{\B}p={\sf c}_k{\sf s}_{\rho^{-1}}^{\B}{\sf s}_{\rho}^{\B}(q_1\odot {\sf c}_k q_2)={\sf c}_k(q_1\odot {\sf c}_kq_2)$$
We next calculate ${\sf c}_k{\sf s}_{\rho^{-1}}p$.
Let $\mu$ be a permutation of $I$ such that $\mu\rho^{-1}I\subseteq I$. Let $j=\mu(\{k\}\cap \rho^{-1}I)$.
Then applying (*), we have:
\begin{equation*}
\begin{split}
&{\sf c}_k{\sf s}_{\rho^{-1}}p={\sf s}^{\B}_{\mu^{-1}}{\sf c}_{j}{\sf s}_{(\mu\rho^{-1}|I)}^{\A}p,\\
&={\sf s}^{\B}_{\mu^{-1}}{\sf c}_{j}{\sf s}_{(\mu\rho^{-1}|I)}^{\A}
{\sf s}^{\A}_{\rho\sigma_1\upharpoonright I}p_1\odot {\sf c}_{(\rho\{k\}\cap \sigma_2I)} {\sf s}^{\B}_{(\rho\sigma_2\upharpoonright I)}p_2,\\
 &={\sf s}^{\B}_{\mu^{-1}}{\sf c}_{j}[{\sf s}_{\mu \sigma_1\upharpoonright I}p_1\odot r].\\
\end{split}
\end{equation*}
where 
$$r={\sf s}_{\mu\rho^{-1}}^{\B}{\sf c}_j {\sf s}_{\rho \sigma_2\upharpoonright I}^{\A}p_2.$$
Now ${\sf c}_kr=r$. Hence, applying the axiom in the small algebra, we get: 
$${\sf s}^{\B}_{\mu^{-1}}{\sf c}_{j}[{\sf s}_{\mu \sigma_1\upharpoonright I}^{\A}p_1]\odot {\sf c}_k q_2
={\sf s}^{\B}_{\mu^{-1}}{\sf c}_{j}[{\sf s}_{\mu \sigma_1\upharpoonright I}^{\A}p_1\odot r].$$
But
$${\sf c}_{\mu(\{k\}\cap \rho^{-1}I)}{\sf s}_{(\mu\sigma_1|I)}^{\A}p_1=
{\sf c}_{\mu(\{k\}\cap \sigma_1I)}{\sf s}_{(\mu\sigma_1|I)}^{\A}p_1.$$
So 
$${\sf s}^{\B}_{\mu^{-1}}{\sf c}_{k}[{\sf s}_{\mu \sigma_1\upharpoonright I}^{\A}p_1]={\sf c}_kq_1,$$ and 
we are done.
To show that neat reducts commute with forming subalgebras, we proceed as in the previous proof replacing finite transformation by transformation.
\end{enumarab}
\end{demo}
We write an $MV$ algebra if the subscripts are clear from context or immaterial.
\begin{definition} Let $\alpha\subseteq \beta$ be infinite sets, $G_{\beta} $ a semigroup of transformations on $\alpha$, $T\subseteq \wp(\alpha)$ 
and $\A$ an  $MV_{G,T}$. An $MV$ algebra $\B$ such that $\A\subseteq \Nr_{\alpha}\B$ is called a $\beta$ dilation of $\A$.
If $A$ generates $\B$ then $\B$ is called a minimal $\beta$ dilation of $\A$.
\end{definition} 
 
The above theorem says that for the  classes we investigate, algebras have a minimal dilation. 
This is an algebraic reflection of adding constants that will act as witnesses to the existential 
quantifiers.
In what follows, until further notice,  we deal with the case when $G$ is the semigroup of all finite transformations on $\alpha$. In this case we have 
$T=\wp_{\omega}(\alpha)$ and we stipulate that $\alpha\sim \Delta x$ is infinite for all $x$ in algebras considered.
To deal with such a case, we need to define certain free algebras, called dimension restricted. The free algebras defined the usual way, will have
the dimensions sets of their elements equal to their dimension, but we do not want that. This concept of dimension restricted free algebras conquers this 
difficulty. It is a deep concept due to Tarski that we borrow from cylindric algebras.
For a class $K$, ${\bf S}$ stands for the operation of forming subalgebras of $K$, and ${\bf P}K$ that of forming direct products.    
\begin{definition}
Let $\delta$ be a cardinal. Let $\alpha$ be an ordinal. Let $T=\wp_{\omega}(\alpha)$ and $G$ be the semigroup of finite transformations on $\alpha$.
Let$_{\alpha} \Fr_{\delta}$ be the absolutely free algebra on $\delta$
generators and of type $MV_{G,T}.$ For an algebra $\A,$ we write
$R\in Con\A$ if $R$ is a congruence relation on $\A.$ Let $\rho\in
{}^{\delta}\wp(\alpha)$. Let $L$ be a class having the same
similarity type as $MV_{G,T}.$ Let
$$Cr_{\delta}^{(\rho)}L=\bigcap\{R: R\in Con_{\alpha}\Fr_{\delta},
{}_{\alpha}\Fr_{\delta}/R\in \mathbf{SP}L, {\mathsf
c}_k^{_{\alpha}\Fr_{\delta}}{\eta}/R=\eta/R \text { for each }$$
$$\eta<\delta \text
{ and each }k\in \alpha\smallsetminus \rho(\eta)\}$$ and
$$\Fr_{\delta}^{\rho}L={}_{\alpha}\Fr_{\beta}/Cr_{\delta}^{(\rho)}L.$$
\end{definition}
The ordinal $\alpha$ does not figure out in $Cr_{\delta}^{(\rho)}L$
and $\Fr_{\delta}^{(\rho)}L$ though it is involved in their
definition. However, $\alpha$ will be clear from context so that no
confusion is likely to ensue.

In what follows $Hom(\A,\B)$ is the set of all homomorphisms from $\A$ to $\B$.

\begin{definition} Assume that $\delta$ is a cardinal, $L\subseteq MV_{G,T}$, $\A\in L$,
$x=\langle x_{\eta}:\eta<\beta\rangle\in {}^{\delta}A$ and $\rho\in
{}^{\delta}\wp(\alpha)$. We say that the sequence $x$ $L$-freely
generates $\A$ under the dimension restricting function $\rho$, or
simply $x$ freely generates $\A$ under $\rho,$ if the following two
conditions hold:
\begin{enumroman}
\item $\A=\Sg^{\A}Rg(x)$ and $\Delta^{\A} x_{\eta}\subseteq \rho(\eta)$ for all $\eta<\delta$.
\item Whenever $\B\in L$, $y=\langle y_{\eta}, \eta<\delta\rangle\in
{}^{\delta}\B$ and $\Delta^{\B}y_{\eta}\subseteq \rho(\eta)$ for
every $\eta<\delta$, then there is a unique $h\in Hom(\A, \B)$ such
that $h\circ x=y$.
\end{enumroman}
\end{definition}
The following theorem can be easily distilled from the literature.

\begin{theorem}
Assume that $\delta$ is a cardinal, $L\subseteq MV_{G,T}$,
$\A\in L$, $x=\langle x_{\eta}:\eta<\delta\rangle\in {}^{\delta}A$
and $\rho\in {}^{\delta}\wp(\alpha).$ Then the following hold:
\end{theorem}
\begin{enumroman}
\item $\Fr_{\delta}^{\rho}L\in MV_{G,T}$
and $x=\langle \eta/Cr_{\delta}^{\rho}L: \eta<\delta \rangle$
$\mathbf{SP}L$- freely generates $\A$ under $\rho$.
\item In order that
$\A\cong \Fr_{\delta}^{\rho}L$ it is necessary and sufficient that
there exists a sequence $x\in {}^{\delta}A$ which $L$ freely
generates $\A$ under $\rho$.
\end{enumroman}
\begin{demo}{Proof} \cite{HMT1} theorems 2.5.35, 2.5.36, 2.5.37.
\end{demo}
We give a metalogical interpretation of such dimension restricted free algebras:
\begin{theorem} Let $\Lambda=(\alpha, P, \rho)$ with $P$ and $\rho$ having common domain $\beta$, 
be a language. Assume that $\alpha\sim \rho(p)$ is infinite for all
$p\in P.$
Then $\Fm/\equiv$ is isomorphic to $\Fr_{\beta}^{\rho}MV_{G,T}.$
Furthermore the isomorphism can be chosen to take atomic formulas to the generators of $\Fr_{\beta}^{\rho}MV_{G,T}$.
\end{theorem}
\begin{demo}{Proof} We have $\Fm/\equiv$ is in $MV_{G,T}$. We have to show that for any  algebra $\A$, 
$a\in {}^{\beta}A$ such that $\Delta a_i\subseteq \rho(i)$, then there exists $k\in Hom(\Fm/\equiv, \A)$ such that 
$h(R_k(v_0,\ldots v_j, \ldots )_{j<\rho(k)}/\equiv)=a_i$ for all $i<\beta$.
Looking at $Fm$ as the universe of the absolutely free algebra  $\Fm$ of $MV_{G,T}$ type, there exists
$h\in Hom(\Fm, \A)$ such that $h(R_k(v_0,\ldots v_j\ldots )_{j<\rho(k)})=a_i$ for each $i<\beta$. 
It suffices to show that $\phi\equiv \psi$ implies that $h(\phi)=h(\psi)$, 
or that $\vdash \phi$ implies $h(\phi)=1$. Let $\Delta =\{\phi\in Fm: h(\phi)=1\}$. 
Then it is tedious but routine to check that $\Delta$ contains all the axioms and is closed under the rules of inference.
\end{demo}
 The following lemma is known:
\begin{lemma} \begin{enumerate}
\item Up to isomorphism, every $MV$ algebra $\A$ is an algebra of $[0,1]^*$ functions over $Spec \A$, where $[0,1]^*$ is an ultrapower of $[0,1]$ 
and  $Spec \A$ is the dual space of $\A$ which is a compact Hausdorff space whose underlying set consists of prime ideals of $\A$.
The ultrapower depends on the cardinality of $A$.
\item An $MV$ algebra is simple if and only if it is isomorphic to a subalgebra of $[0,1]$. 
An $MV$ algebra is semisimple iff it is isomorphic to a separating $MV$ algebras of $[0,1]$ 
valued functions on some compact Hausdorff space, with pointwise operations.
\end{enumerate}
\end{lemma}\begin{definition} Let $\A$ be an $MV$ algebra. A filter in $\A$ is a subset of $F$ of $\A$ such that
\begin{enumroman}
\item if $a,b\in F$, then $a\odot b\in F,$
\item if $a\in F$ and $b\geq a$, then $b\in F$.
\end{enumroman}
\end{definition}

The filter $F$ generated by $X$, is the set of elements $\{x\in A: \exists a_0, a_1\ldots a_n\in X: x\geq a_0\odot a_1\ldots a_n\}$.
We write $\prod_{i<n}a_i$ for $a_0\odot\ldots a_{n-1}.$

\section{Interpolation Theorems}  

In this section, we prove algebraically our interpolation theorems.
The following definition is the algebraic counterpart of the interpolation property:
\begin{definition}
An algebra $\A$ has the interpolation property if for all $X_1, X_2\subseteq \A$ if $a\in \Sg^{\A}X_1$ and $b\in \Sg^{\A}X_2$ are such that $a\leq b$, 
then there exists $c\in \Sg^{\A}(X_1\cap X_2)$ such that $a\leq c\leq b$.
\end{definition}
Now we are ready for:
\begin{theorem}\label{in} Let $\alpha$ be an infinite set. Let $\mu$ be a cardinal.  Let $\rho:\mu\to \wp(\alpha)$ such that
$\alpha\sim \rho(i)$ is infinite for all $i\in \mu$. Let $G$ be the semigroup of finite transformations on $\alpha$, and let $T=\wp_{\omega}(\alpha)$.
Then $\Fr_{\beta}^{\rho}MV_{G,T}$ has the interpolation property.
\end{theorem}
\begin{enumerate}
\item We assume that $\alpha$ is an ordinal. We first show that for any ordinal $\beta>\alpha$,
the sequence $\langle \eta/Cr_{\mu}^{\rho}MV_{G,T}:
\eta<\mu\rangle$ $MV_{G,T}$ - freely generates
$\Nr_{\alpha}\Fr_{\mu}^{\rho}(MV_{G,\bar{T}})$ where $\bar{T}=\wp_{\omega}(\beta).$  Let $\B\in
MV_{G,T}$ and $a=\langle a_{\eta} :\eta<\mu\rangle\in {}^{\mu}B$
be such that $\Delta a_{\eta}\subseteq \rho(\eta)$ for all
$\eta<\mu$.  In  $\Fr_{\mu}^{\rho}(MV_{G,\bar{T}})$ we have $\alpha\sim \rho x$ is infinite. 
Assuming that $Rg a$ generates $\B$, we have for all $x\in \B$, $\Delta x\sim \rho(x)$ is infinite.
Therefore, by theorem \ref{net},  $\B$ neatly embeds in an algebra in $MV_{\beta, \bar{T}}$. Let
$\D=\Fr_{\mu}^{\rho}(MV_{G,\bar{T}})$. Then, we claim that $x=\langle
\eta/Cr_{\mu}^{\rho}MV_{G,\bar{T}}:\eta<\mu\rangle\in {}^{\mu}D,$
$S\Nr_{\alpha}MV_{G\,\bar{T}}$  freely generates
$\Sg^{\Nr_{\alpha}\D}Rgx$. Indeed, consider $\C\in
\Nr_{\alpha}MV_{G,\bar{T}}$ and $y\in {}^{\mu}C$ such that $\Delta
y_{\eta}\subseteq \rho \eta$ for all $\eta<\mu$. Let $\C'\in
MV_{G,\bar{T}}$ be such that $\C=\Nr_{\alpha}\C'$. Then clearly $y\in
{}^{\mu}C'$ and $\Delta y_{\eta}\subseteq \alpha$ for all
$\eta<\mu$. Then there exists $h\in Hom(\D, \C')$ such that $h\circ
x=y$. Hence $h\in Hom(\Rd_{\alpha}\D, \Rd_{\alpha}\C')$, thus $h\in
Hom(\Sg^{\Rd_{\alpha}\D}Rg x, \Sg^{\Rd_{\alpha}\C'}h(Rgx))$. Since
$Rg x\subseteq Nr_{\alpha} D$, we have $h\in
Hom(\Sg^{\Nr_{\alpha}\D}Rgx, \C)$. The conclusion now follows. Therefore
there exists $h:\Sg^{\Nr_{\alpha}\D}Rgx\to \B$ such that
$h(\eta/Cr_{\mu}^{\rho}MV_{G,T}{\beta})= a_{\eta}$. But by theorem \ref{net}, we have 
$$\Nr_{\alpha}\Fr_{\mu}^{\rho}(MV_{G,\bar{T}})=\Nr_{\alpha}(\Sg^{\D}Rgx)=\Sg^{\Nr_{\alpha}\D}Rgx.$$
Therefore, as claimed,  $\langle \eta/Cr_{\mu}^{\rho}(MV_{G,\bar{T}}):
\eta<\mu\rangle$ $MV_{G,T}$ freely generates
$\Nr_{\alpha}\Fr_{\mu}^{\rho}(MV_{G,\bar{T}})$, so that 
$$\Nr_{\alpha}\Fr_{\mu}^{\rho}(MV_{G_{\beta},\bar{T}})\cong \Fr_{\mu}^{\rho}(MV_{G_{\alpha},T}).$$

\item The idea we implement here is basically the one used in \cite{IGPL}, but the details are far from being identical. 
We have to check that the proof goes through in the absence of distributivity and idempotency 
enjoyed by  Boolean algebras. Surprisingly it does, but with non-trivial modifications.
Let $\B=\Fr_{\beta}^{\rho}(MV_{G,T})$. Let $a\in \Sg^{\B}X_1$ and $b\in \Sg^{\B}X_2$ be such that $a\leq b$. 
Notice that this is the lattice order, and that $a\leq b$ iff $a\odot \neg b=0.$
We want to find an interpolant in 
$\Sg^{\B}(X_1\cap X_2)$. Assume that $\kappa$ is a regular cardinal $>max (|A|,|\alpha|)$.   
Let $\C=\Fr_{\beta}^{\rho}MV_{G,\bar{T}}$, where $\bar{T}=\wp_{\omega}(\kappa)$. Then by what is proven in the previous item, we have
$\B=\Nr_{\alpha}\C$, 
and $B$ generates $\C$.
If an interpolant exists in the big algebra $\C$, then an interpolant exists in the smaller one $\B$. 
For assume there exists $c\in \Sg^{\C}(X_1\cap X_2)$ 
such that that $a\leq c\leq b.$ Then there exists a finite 
$\Gamma\subseteq \kappa\sim \alpha$ such that $a\leq {\sf c}_{(\Gamma)}c\leq b$. Now 
by theorem \ref{net} we have 
$${\sf c}_{(\Gamma)}c\in \Nr_{\alpha}\Sg^{\C}(X_1\cap X_2)=S\g^{\Nr_{\alpha}\C}(X_1\cap X_2)=\Sg^{\B}(X_1\cap X_2).$$
So assume that no interpolant exists in $\B$, then no interpolant exists in $\C$. We will reach a contradiction. Recall that 
$\kappa>max(|A|, \alpha).$
Arrange $\kappa\times \Sg^{\C}(X_1)$ 
and $\kappa\times \Sg^{\C}(X_2)$ 
into $\kappa$-termed sequences
$$\langle (k_i,x_i): i\in \kappa\rangle\text {  and  }\langle (l_i,y_i):i\in \kappa\rangle
\text {  respectively.}$$ Since $\kappa$ is regular, we can define by recursion 
$\kappa$-termed sequences 
$$\langle u_i:i\in \kappa\rangle \text { and }\langle v_i:i\in \kappa\rangle,$$ 
such that for all $i\in \kappa$ we have:
$$u_i\in \kappa\sim
(\Delta a\cup \Delta b)\cup \cup_{j\leq i}(\Delta x_j\cup \Delta y_j)\cup \{u_j:j<i\}\cup \{v_j:j<i\}$$
and
$$v_i\in \kappa\sim\Delta a\cup \Delta b)\cup 
\cup_{j\leq i}(\Delta x_j\cup \Delta y_j)\cup \{u_j:j\leq i\}\cup \{v_j:j<i\}.$$
For an $MV$ algebra $\D$  and $Y\subseteq \D$, we write 
$fl^{\D}Y$ to denote the $MV$ filter generated by $Y$ in $\C.$ For an $MV_{G,T}$ $\D$, we write $\Rd_{MV}\D$ for its $MV$ reduct obtained by 
discarding the operations of cylindrifiers and substitutions. Now let 
$$Y_1= \{a\}\cup \{-{\sf  c}_{k_i}x_i\oplus{\sf s}_{u_i}^{k_i}x_i: i\in \kappa\},$$
$$Y_2=\{-b\}\cup \{-{\sf  c}_{l_i}y_i\oplus {\sf s}_{v_i}^{l_i}y_i:i\in \kappa\},$$
$$H_1= fl^{\Rd_{MV}\Sg^{\C}(X_1)}Y_1,\  H_2=fl^{\Rd_{MV}\Sg^{\C}(X_2)}Y_2,$$ and 
$$H=fl^{\Rd_{MV}\Sg^{\C}(X_1\cap X_2)}[(H_1\cap \Sg^{\C}(X_1\cap X_2)
\cup (H_2\cap \Sg^{\C}(X_1\cap X_2)].$$
We claim that $H$ is a proper filter of $\Sg^{\C}(X_1\cap X_2).$
To prove this it is sufficient to consider any pair of finite, strictly
increasing sequences of natural numbers 
$$\eta(0)<\eta(1)\cdots <\eta(n-1)<\omega\text { and } \xi(0)<\xi(1)<\cdots 
<\xi(m-1)<\omega,$$
and to prove that the following condition holds:

$(1)$  For any $b_0$, $b_1\in \Sg^{\C}(X_1\cap X_2)$ such that
$$a\odot\prod_{i<n}[a^{l-1}\odot (-{\sf  c}_{k_{\eta(i)}}x_{\eta(i)}\oplus{\sf s}_{u_{\eta(i)}}^{k_{\eta(i)}}x_{\eta(i)})^{l_i}]\leq b_0$$  
and
$$(-b)\odot\prod_{i<m}
[(-b)^{k-1}\odot (-{\sf  c}_{l_{\xi(i)}}y_{\xi(i)}\oplus{\sf s}_{v_{\xi(i)}}^{l_{\xi(i)}}y_{\xi(i)})^{k_i}]\leq b_1,$$
where for $i<n$ and $j<m$, $l_i, k_j$ as well as $l$ and $k$ are finite ordinals $>0$, we have 
$$b_0\odot b_1\neq 0.$$

We prove this by induction on $l-1+n+m+k-1=l+n+m+k-2$, which we can assume is $\geq 0$; for else there would be nothing to 
prove. 
If this number is equal to $0$, then $(1)$ simply 
expresses the fact that no interpolant of $a$ and $b$ exists in 
$\Sg^{\C}(X_1\cap X_2).$
In more detail: if $l+n+m+k-2=0$, then $a\leq b_0$
and $-b\leq b_1$. So if $b_0\odot b_1=0$, we get $a\leq b_0\leq -b_1\leq b$ and $b_0$ would be the desired interpolant.

Now assume that $l+n+m+k-2>0$ and for the time being suppose that $\eta(n-1)>\xi(m-1)$.
Apply ${\sf  c}_{u_{\eta(n-1)}}$ to both sides of the first inclusion of (1).
By $u_{\eta(n-1)}\notin \Delta a$, i.e. ${\sf  c}_{u_{\eta(n-1)}}a=a$, 
and by noting that ${\sf  c}_i({\sf  c}_ix\odot y)={\sf  c}_ix\odot {\sf  c}_iy$, we get (2)
$$a\odot {\sf  c}_{u_{\eta(n-1)}}\prod_{i<n} [a^{l-1}\odot ( -{\sf  c}_{k_{\eta(i)}}x_{\eta(i)}
\oplus{\sf s}_{u_{\eta(i)}}^{k_{\eta(i)}}x_{\eta(i)})^{l_i}]\leq {\sf  c}_{u_{\eta(n-1)}}b_0.$$
 Now apply ${\sf  q}_{u_{\eta(n-1)}}$ to the second inclusion of (1).
By observing that $u_{\eta(n-1)}\notin \Delta b=\Delta (-b)$ we get (3)
$$(-b)\odot \prod_{j<m}[(-b)^{k-1}\odot ({\sf  q}_{u_{\eta(n-1)}}(-{\sf  c}_{l_{\xi(j)}}y_{\xi(j)}\oplus{\sf s}_{v_{\xi(j)}}^{l_{\xi(j)}}
y_{\xi(j)})^{k_i}]\leq {\sf  q}_{u_{\eta(n-1)}}b_1.$$
Before going on, we formulate (and prove) a claim that will enable us to eliminate
the quantifier ${\sf  c}_{u_{\eta(n-1)}}$ (and its dual) from (2) (and (3)) above.

\begin{athm}{Claim 1} Let $t_j=-{\sf  c}_{l_{\xi(j)}}y_{\xi(j)}\oplus{\sf s}_{v_{\xi(j)}}^{l_{\xi(j)}}y_{\xi(j)}.$ Then 
$${\sf  q}_{u_{\eta(n-1)}}t_j=t_j\text { for all }j<m.$$
\end{athm}
\begin{demo}{Proof of Claim 1}
Let $j<m$ . Then we have 
\begin{equation*}
\begin{split}
{\sf  q}_{u_{\eta(n-1)}}(-{\sf  c}_{l_{\xi(j)}}y_{\xi(j)})
&=-{\sf  c}_{l_{\xi(j)}}y_{\xi(j)}\\ 
{\sf  c}_{u_{\eta(n-1)}}^{\partial}
({\sf s}_{v_{\xi(j)}}^{l_{\xi(j)}}y_{\xi(j)})
&={\sf s}_{v_{\xi(j)}}^{l_{\xi(j)}}y_{\xi(j)}.\\
\end{split}
\end{equation*}
Indeed,  computing we get
\begin{equation*}
\begin{split}
{\sf  q}_{u_{\eta(n-1)}}(-{\sf  c}_{l_{\xi(j)}}y_{\xi(j)})
&=-{\sf  c}_{u_{\eta_{n-1}}}-(-{\sf  c}_{l_{\xi(j)}}y_{\xi(j)})\\
&= -{\sf  c}_{u_{\eta(n-1)}}{\sf  c}_{l_{\xi(j)}}y_{\xi(j)}\\
&=-{\sf  c}_{l_{\xi(j)}}y_{\xi(j)}.
\end{split}
\end{equation*}
Similarly,  we have
\begin{equation*}
\begin{split}
{\sf  q}_{u_{\eta(n-1)}} ({\sf s}_{v_{\xi(j)}}^{l_{\xi(j)}}y_{\xi(j)})
&=-{\sf  c}_{u_{\eta(n-1)}}- ({\sf s}_{v_{\xi(j)}}^{l_{\xi(j)}}y_{\xi(j)})\\
&=-{\sf  c}_{u_{\eta(n-1)}} ({\sf s}_{v_{\xi(j)}}^{l_{\xi(j)}}-y_{\xi(j)})\\
&=- {\sf s}_{v_{\xi(j)}}^{l_{\xi(j)}}-y_{\xi(j)}\\
&= {\sf s}_{v_{\xi(j)}}^{l_{\xi(i)}}y_{\xi(j)}.
\end{split}
\end{equation*}
By ${\sf  q}_i({\sf  q}_i x\oplus y)= {\sf q}_ix\oplus{\sf  q}_i y$
we get from the above that
\begin{equation*}
\begin{split} 
{\sf  q}_{u_{\eta(n-1)}}(t_j)&={\sf  q}_{u_{\eta(n-1)}}(-{\sf  c}_{l_{\xi(j)}}y_{\xi(j)}\oplus{\sf s}_{v_{\xi(j)}}^{l_{\xi(j)}}y_{\xi(j)})\\
&={\sf  q}_{u_{\eta(n-1)}}-{\sf  c}_{l_{\xi(j)}}y_{\xi(j)}\oplus{\sf  q}_{u_{\eta(n-1)}}
{\sf s}_{v_{\xi(j)}}^{l_{\xi(j)}}y_{\xi(j)}\\
&=-{\sf  c}_{l_{\xi(j)}}y_{\xi(j)}\oplus {\sf s}_{v_{\xi(j)}}^{l_{\xi(j)}}y_{\xi(j)}=t_j.
\end{split}
\end{equation*}
 
\end{demo}
\begin{athm}{Claim 2} For  each $i<n$ and each $j<m$, let
$$z_i=-{\sf  c}_{k_{\eta(i)}}x_{\eta(i)}\oplus{\sf s}_{u_{\eta(i)}}^{k_{\eta(i)}}x_{\eta(i)}.$$  Then
$${\sf  c}_{u_{\eta(n-1)}}z_i=z_i\text { for }i<n-1 \text { and }{\sf  c}_{u_{\eta(n-1)}}z_{n-1}=1.$$
\end{athm}
\begin{demo}{Proof of Claim 2}
Let $i<n-1$. Then by the choice of witnesses we have 
$$u_{\eta(n-1)}\neq u_{\eta(i)}.$$ 
Also it is easy to see that for all $i,j\in \alpha$ we have 
$$\Delta {\sf  c}_jx\subseteq \Delta x\text {  and that }
\Delta {\sf s}_j^ix\subseteq \Delta x\smallsetminus\{i\}\cup \{j\},$$
In particular, 
$$u_{\eta(n-1)}\notin \Delta {\sf  c}_{k_{\eta(i)}}x_{\eta(i)}\text { and }
u_{\eta(n-1)}\notin \Delta ({\sf s}_{u_{\eta(i)}}^{k_{\eta(i)}}x_{\eta(i)}).$$
It thus follows that 
$${\sf  c}_{u_{\eta(n-1)}}(-{\sf  c}_{k_{\eta(i)}}x_{\eta(i)})=-{\sf  c}_{k_{\eta(i)}}x_{\eta(i)}\text { and }
{\sf  c}_{u_{\eta(n-1)}} ({\sf s}_{u_{\eta(i)}}^{k_{\eta(i)}}x_{\eta(i)})={\sf s}_{u_{\eta(i)}}^{k_{\eta(i)}}
x_{\eta(i)}.$$
Finally, by properties of ${\sf  c}_{u_{\eta(n-1)}}$, we get 
$${\sf  c}_{u_{\eta(n-1)}} z_i=z_i \text { for }  i<n-1.$$

{\it Proof of ${\sf  c}_{u_{\eta(n-1)}}z_{n-1}=1.$}

Computing we get, by $u_{\eta(n-1)}\notin \Delta x_{\xi(n-1)},$ 
and by familiar axioms of substitutions, namely items 7,8,9 in theorem \ref{d},
the following:
\begin{equation*}
\begin{split}
&{\sf  c}_{u_{\eta(n-1)}}(-{\sf  c}_{k_{\eta(n-1)}}x_{\eta(n-1)}
\oplus {\sf s}_{u_{\eta(n-1)}}^{k_{\eta(n-1)}}x_{\eta(n-1)})\\
&={\sf  c}_{u_{\eta(n-1)}}-{\sf  c}_{k_{\eta(n-1)}}x_{\eta(n-1)}\oplus {\sf  c}_{u_{\eta(n-1)}}
{\sf s}_{u_{\eta(n-1)}}^{k_{\eta(n-1)}}x_{\eta(n-1)}\\
&=-{\sf  c}_{k_{\eta(n-1)}}x_{\eta(n-1)}\oplus {\sf  c}_{u_{\eta(n-1)}}{\sf s}_{u_{\eta(n-1)}}^{k_{\eta(n-1)}}
x_{\eta(n-1)}\\
&=-{\sf  c}_{k_{\eta(n-1)}}x_{\eta(n-1)}\oplus {\sf  c}_{u_{\eta(n-1)}}{\sf s}_{u_{\eta(n-1)}}^{k_{\eta(n-1)}}
{\sf  c}_{u_{\eta(n-1)}}x_{\eta(n-1)}\\
&=-{\sf  c}_{k_{\eta(n-1)}}x_{\eta(n-1)}\oplus {\sf  c}_{k_{\eta(n-1)}}{\sf s}_{k_{\eta(n-1)}}^{u_{\eta(n-1)}}
{\sf  c}_{u_{\eta(n-1)}}x_{\eta(n-1)}\\
&=-{\sf  c}_{k_{\eta(n-1)}}x_{\eta(n-1)}\oplus {\sf  c}_{k_{\eta(n-1)}}{\sf  c}_{u_{\eta(n-1)}}x_{\eta(n-1)}\\
&= -{\sf  c}_{k_{\eta(n-1)}}x_{\eta(n-1)}\oplus {\sf  c}_{k_{\eta(n-1)}}x_{\eta(n-1)}=1.\\
\end{split}
\end{equation*}
\end{demo}

By the above proven claims, and using the notation introduced in claim 2,   we have 
\begin{equation*}
\begin{split}
&{\sf  c}_{u_{\eta(n-1)}}(a^{l-1}\odot \prod_{i<n}z_i^{l_i})\\
&=a^{l-1}\odot {\sf  c}_{u_{\eta(n-1)}}\prod_{i<n-1}(z_i^{l_i}\odot z_{n-1}^{l_{n-1}})\\
&=a^{l-1}\odot {\sf  c}_{u_{\eta(n-1)}}\prod_{i<n-1}z_i^{l_i}\odot
({\sf  c}_{u_{\eta(n-1)}}z_{n-1}^{l_{n-1}})\\
&=a^{l-1}\odot {\sf  c}_{u_{\eta(n-1)}}\prod_{i<n-1}z_i^{l_i}\odot 
({\sf  c}_{u_{\eta(n-1)}}z_{n-1})^{l_{n-1}}\\
&=a^{l-1}\odot \prod_{i<n-1}z_i^{l_i}.\\
\end{split}
\end{equation*}
Combined with (2) we obtain 
$$a\odot \prod_{i<n-1}[a^{l-1}\odot ( -{\sf  c}_{k_{\eta(i)}}x_{\eta(i)}\oplus{\sf s}_{u_{\eta(i)}}^{k_{\eta(i)}}x_{\eta(i)})^{l_i}]
 \leq {\sf  c}_{u_{\eta(n-1)}}b_0.$$
On the other hand, from our proven claims  and (3), 
 it follows that
$$(-b)\odot \prod_{j<m}
[(-b)^{k-1}\odot (-{\sf  c}_{l_{\xi(j)}}y_{\xi(j)}\oplus{\sf s}_{v_{\xi(j)}}^{l_{\xi(j)}}y_{\xi(j)})^{k_i}]\leq {\sf  q}_{u_{\eta(n-1)}}b_1.$$
Now making use of the induction hypothesis, we get
$${\sf  c}_{u_{\eta(n-1)}}b_0\odot {\sf q}_{u_{\eta(n-1)}}b_1\neq 0;$$ 
and hence that
$$b_0\odot {\sf  q}_{u_{\eta(n-1)}}b_1\neq 0.$$ 
From 
$$b_0\odot {\sf  q}_{u_{\eta(n-1)}}b_1\leq b_0\odot b_1,$$
we reach the desired conclusion, i.e. that 
$$b_0\odot b_1\neq 0.$$
The other case, when $\eta(n-1)\leq \xi(m-1)$ can be treated analgously
and is therefore left to the reader.
We have proved that $H$ is a proper filter.

\item Proving that $H$ is a proper filter of $\Sg^{\C}(X_1\cap X_2)$, 
let $H^*$ be a maximal filter of $\Sg^{\C}(X_1\cap X_2)$ 
containing $H.$ 
We obtain  maximal $F_1$ and $F_2$ of $\Sg^{\C}(X_1)$ and $\Sg^{\C}(X_2)$, 
respectively, such that 
$$H^*\subseteq F_1,\ \  H^*\subseteq F_2$$
and (**)
$$F_1\cap \Sg^{\C}(X_1\cap X_2)= H^*= F_2\cap \Sg^{\C}(X_1\cap X_2).$$
Now for all $x\in \Sg^{\C}(X_1\cap X_2)$ we have 
$$x\in F_1\Longleftrightarrow x\in F_2.$$ 
Also from how we defined our maximal  filters, $F_i$ for $i\in \{1,2\}$ satisfy the following
condition:

(*) For all $k<\kappa$, for all $x\in \Sg^{\C}X_i$ 
if ${\sf  c}_kx\in F_i$ then ${\sf s}_l^kx$ is in $F_i$ for some $l\notin \Delta x.$  

Let $V$ be the set of all finite transformations on $\kappa$, i.e
$$V=\{\sigma\in {}^{\kappa}\kappa: |\{\sigma(i)\neq i\}|<\omega\}.$$
Every $\tau\in V$ defines a unary operation ${\sf s}_{\tau}$ on $\C$ as follows (see the proof of theorem \ref{net}).
If $\tau=[u_0|v_0, u_1|v_1,\ldots,
u_{k-1}|v_{k-1}]$ is a finite transformation, if $x\in C$ and if
$\pi_0,\ldots ,\pi_{k-1}$ are in this order the first $k$ ordinals
in $\kappa\sim (\Delta x\cup Rg(u)\cup Rg(v))$, then
$${\mathsf s}_{\tau}x={\mathsf s}_{v_0}^{\pi_0}\ldots
{\mathsf s}_{v_{k-1}}^{\pi_{k-1}}{\mathsf s}_{\pi_0}^{u_0}\ldots
{\mathsf s}_{\pi_{k-1}}^{u_{k-1}}x.$$
Let ${\D}_i=\Sg^{\C}X_i$, $i=1,2.$
Then let  $\psi_i$ be the map defined as follows:
$$\psi_i: \D_i\to \F(V, \D_i/F_i)$$
$$\psi_i(a)(x)={\sf s}_x^{\C}a/F_i.$$  
For brevity, we omit the superscript $i$, so that 
$$\psi: \D\to \F(V, \D/F).$$
We have $\D/F$ is a simple algebra, hence it is isomorphic to a subalgebra of $[0,1]$.
So we consider $\psi$ as a mapping from $\D$ into $\F(V,[0,1]).$
We check that it is a homomorphism.
We only check cylindrifiers (the other operations are straightforward to check.) Now, we have
$$\psi({\sf c}_ka)(x)={\sf s}_x{\sf c}_ka/F.$$
Let $$l\in \{\mu\in \kappa: x^{-1}\{\mu\}=\{\mu\}\}\sim \Delta a.$$
Such an $l$ clearly exists.
Let $$\tau=x\circ [k,l].$$
Then by familiar substitution rules we have
$${\sf c}_l{\sf s}_{\tau}a={\sf s}_{\tau}{\sf c}_ka={\sf s}_x{\sf c}_ka,$$
and by the choice of $F,$ we have
$${\sf c}_l{\sf s}_{\tau}a\in F\Longleftrightarrow {\sf s}^l_u{\sf s}_{\tau}a\in F.$$
We use the following helpful notation. 
For a function $f$, the function $g=f(a\to u)$ is defined by $g(x)=f(x)$ for $x\neq a$ and $g(a)=u$. Now we have
\begin{equation*}
\begin{split}
\psi({\sf c}_ka)(x)
&={\sf s}_x{\sf c}_ka/F\\
&={\sf c}_l{\sf s}_{\tau}a/F\\
&={\sf s}_u^l{\sf s}_{\tau}a/F \\
&={\sf s}_{x(k\to u)}a/F\\
&=\psi a({x(k\to u}))/F\\
&\leq {\sf c}_k\psi(a)(x)
\end{split}
\end{equation*}

Conversely if $y\equiv_k x$, then 
$$\psi(a)(y)={\sf s}_ya/F\leq {\sf s}_x{\sf c}_ka/F=\psi({\sf c}_ka)(x).$$
We have proved that $\psi$ is a homomorphism. Then putting back superscripts, we show that $\psi_1$ and $\psi_2$ agree on their common part
$\Sg^{\C}(X_1\cap X_2)$.

Le $a\in \Sg^{\C}(X_1\cap X_2).$ Then
\begin{equation*}
\begin{split}
\psi_1(a)(x)
&={\sf s}_xa/F_1\\
&={\sf s}_xa/F_1\cap \Sg^{\C}(X_1\cap X_2)\\
&={\sf s}_xa/H^*\\
&={\sf s}_xa/F_2\cap \Sg^{\C}(X_1\cap X_2)\\
&={\sf s}_xa/H^*\\
&=\psi_2(a)(x).
\end{split}
\end{equation*}
Assuming that $X_1\cup X_2$ generates $\C$, we  have $\psi_1\cup \psi_2$ defines a function on $\C$ into $\F(V,[0,1])$, 
since they agree on the common part; 
by freeness they can be pasted to give $\psi:\C\to \F(V,[0,1])$ 
such that $\psi(a\odot -b)\neq 0$, because the identity substitution is in $\psi(a\odot -b)$ by definition, 
but this contradicts that $a\leq b$ and the proof is complete.

\end{enumerate}
Reformulating definition \ref{rep}, we have:
\begin{definition} An $MV_{G,T}$ algebra is representable if it is isomorphic to a subdirect product of algebras of the form $\F(V, [0,1])$ where
$V\subseteq {}^IX$, for some sets $I$ and $X$.
\end{definition}
\begin{corollary} Let $\alpha$ be infinite, $G$ be the semigroup of finite transformations on $\alpha$ and $T=\wp_{\omega}(\alpha).$ Then any 
$\A$ in $MV_{G,T}$, such that $\alpha\sim \Delta x$, is representable.
\end{corollary}
\begin{demo}{Proof} Let $\A$ be given and $a\neq 0$ be in $A$. Let $\kappa$ be a regular cardinal $>max(|\alpha|, |A|)$.  
Let $\B\in MV_{G_{\kappa},T}$ such that $\A=\Nr_{\alpha}\B$. Let $\langle (k_i,x_i): i\in \kappa \rangle$ be an enumeration of $\kappa\times B.$
Since  $\kappa$ is regular, we can define by recursion a
$\kappa$-termed sequence  $\langle u_i:i\in \kappa\rangle$ 
such that for all $i\in \kappa$ we have:
$u_i\in \kappa\sim
(\Delta a\cup \bigcup_{j\leq i}\Delta x_j\cup \{u_j:j<i\}).$
Let  $Y=\{a\}\cup \{-{\sf  c}_{k_i}x_i\oplus{\sf s}_{u_i}^{k_i}x_i: i\in \kappa\}.$ Let $H$ be the filter generated by $Y$; 
then $H$  is proper, take the maximal filter containing $H$ and $a$,
and define $\psi(b)x={\sf s}_xb/F$ where $b\in B$ and $x\in V$,  and $V$ is as defined in the previous proof.
Then $\psi(a)\neq 0$, and $\psi$ establishes the representability of $\B$, hence of $\A$.
\end{demo}

\begin{theorem} Let $\alpha$ be an infinite set.  \begin{enumarab}
\item If $\alpha$ is countable, $G$ is a rich semigroup on $\alpha$ and $T=\wp_{\omega}(\alpha)$ 
then $\Fr_{\beta}MV_{G,T}$ has the interpolation property.
\item  If $\alpha$ is arbitrary, $G={}^{\alpha}\alpha$ and
$T=\wp_{\omega}(\alpha)$, then $\Fr_{\beta}MV_{G,T}$ has the interpolation property.
\end{enumarab}
\end{theorem}
\begin{demo}{Proof} 
\begin{enumarab}
\item The proof is similar to the proof of theorem \ref{in}. We can assume that $\alpha=\omega$. Let $\A=\Fr_{\beta}MV_{G,T}$. By theorem \ref{rich}
let $\B\in MV_{\bar{G},\bar{T}}$ such that $\A=\Nr_{\omega}\B$. Though irrelevant to the present proof it can be proved using the above reasoning that  
$\B=\Fr_{\beta}MV_{\bar{G}, \bar{T}}$.
Let $X_1,X_2\subseteq A$, assume that $a\in \Sg^{\A}X_1$ and $b\in \Sg^{\A}X_2$ such that $a\leq b$. 
Assume that no interpolant exists in $\A$. Then, as in the proof of theorem \ref{in},  no interpolant exists in $\B$.
Let $\alpha=\omega+\omega$. We can define by recursion
$\omega$-termed sequences of witnesses: 
$$\langle u_i:i\in \omega\rangle \text { and }\langle v_i:i\in \omega\rangle$$ 
such that for all $i\in \omega$ we have:
$$u_i\in \alpha\smallsetminus(\Delta a\cup \Delta b)\cup \cup_{j\leq i}(\Delta x_j\cup \Delta y_j)\cup \{u_j:j<i\}\cup \{v_j:j<i\}$$
and
$$v_i\in \alpha\smallsetminus(\Delta a\cup \Delta b)\cup 
\cup_{j\leq i}(\Delta x_j\cup \Delta y_j)\cup \{u_j:j\leq i\}\cup \{v_j:j<i\}.$$
$$Y_1= \{a\}\cup \{-{\sf c}_{k_i}x_i\oplus {\sf s}_{u_i}^{k_i}x_i: i\in \omega\},$$
$$Y_2=\{-b\}\cup \{-{\sf c}_{l_i}y_i\oplus {\sf s}_{v_i}^{l_i}y_i:i\in \omega\},$$
$$H_1= fl^{\Rd_{MV}\B}(X_1)Y_1,\  H_2=fl^{\Rd_{MV}\B}(X_2)Y_2,$$ and 
$$H=fl^{\Rd_{MV}\B(X_1\cap X_2)}[(H_1\cap \Rd_{MV}\B(X_1\cap X_2)\cup (H_2\cap A^+(X_1\cap X_2)].$$
Then $H$ is a proper filter of $\Sg^{\B}(X_1\cap X_2)$; same reasoning as in the proof of theorem \ref{in}.
Let $V=\bigcup_{\tau\in G}{}^{\omega}\alpha^{\tau}$. Then for $s \in V$, let $\bar{s}= s\cup Id \in \bar{G}$. Then define
for $x\in V$, $\psi_1(a)x={\sf s}_{\bar{x}}^{\B}/F_1$ and $\psi_2={\sf s}_{\bar{x}}^{\B}/F_2$. Then these are homomorphisms that 
can be pasted together to give a homomorphism $\psi$ with domain $\A$ 
such that $\psi(a\odot -b)\neq 0$. This contradicts $a\leq b$.

\item When $G=\wp({\alpha})$, then the first part of the proof is identical to that of theorem \ref{in} but resorting to \ref{net} (2) instead. 
Let $\A=\Fr_{\beta}MV_{G,T}$. 
Let $\kappa$ be a regular cardinal $>max(|A|,\alpha)$, $\B$ be an algebra of dimension $\kappa$ such that $\A=\Nr_{\alpha}\B$
and then proceed as above constructing the proper filter $H$. In defining the representability functions, instead of 
$V$ one takes $^{\kappa}\kappa$ because all substitutions are at hand, and then the proof is exactly the same. 
\end{enumarab}
\end{demo}

\begin{corollary} Let $\alpha$ be an infinite set.  If $\alpha$ is countable, $G$ is a rich semigroup on 
$\alpha$ and $T=\wp_{\omega}(\alpha),$ 
or $\alpha$ is arbitrary, $G={}^{\alpha}\alpha$ and
$T=\wp_{\omega}(\alpha),$ and in both cases $\A\in MV_{G,T},$ then $\A$ is representable.
\end{corollary}
\begin{demo}{Proof} From the above arguments.
\end{demo}
To prove our result for full polyadic algebras (where cylindrifiers are taken on all subsets and all substitutions are available), 
we closely follow $\cite{MLQ}$.
\begin{definition} Let $\A\in MV_{G,T}$ where $G={}^{\alpha}\alpha$ and $T=\wp(\alpha)$. 
\begin{enumroman}

\item If $J\subseteq \alpha$, an element $a\in A$ is independent of $J$ if ${\sf c}_{(J)}p=p$.
$J$ supports $a$ if $a$ is independent of $\alpha\sim J$.
$A_J=\{a\in A: \text { $J$ supports $a$ }\}.$

\item The effective degree of $\A$ is the smallest cardinal $\mathfrak{e}$ such that each element of $\A$ admits a support 
whose cardinality does not exceed $\mathfrak{e}$.
\item The local degree of $\A$ is the smallest cardinal $\mathfrak{m}$ such that each element of $\A$ has cardinality $<\mathfrak{m}$. 
\item The effective cardinality of $\A$ is $\mathfrak{c}=|A_J|$ where $|J|=\mathfrak{e}.$ (This is independent of $J$). 
\end{enumroman}
\end{definition}

\begin{theorem} Let $\alpha$ be infinite, $G={}^{\alpha}\alpha$ and $T=\wp(\alpha)$. Then $\Fr_{\beta}MV_{G,T}$ has the interpolation property.
\end{theorem} 
\begin{demo}{Proof} 
\begin{enumarab}
\item The first part of the proof is identical to that in \cite{MLQ}, but we include it for the sake of completeness referring to op cit for detailed arguments. 
Let $\mathfrak{m}$ be the local degree of $\A=\Fr_{\beta}MV_{G,T}$, $\mathfrak{c}$ its effective cardinality and $\mathfrak{n}$ 
be any cardinal such that $\mathfrak{n}\geq \mathfrak{c}$ 
and $\sum_{s<m}\mathfrak{n}^s=\mathfrak{n}$. Let $X_1, X_2\subseteq A$ $a\in \Sg^{\A}X_1$ and $b\in \Sg^{\A}X_2$ such that $a\leq b$. 
We want to find an interpolant.
By theorem \ref{net},  
there exists $\B\in MV_{G_{\mathfrak{n}},\bar{T}}$  where $\bar{T}=\wp(\kappa)$ such that $\A\subseteq \Nr_{\alpha}\B$ and $A$ generates $\B$. 
(It can be proved that $\B=\Fr_{\beta}MV_{G_{\mathfrak{n},\bar{T}}}).$
Being a minimal dilation of $\A$, the local degree of $\B$ is the same as that of $\A$, in particular each $x\in \B$ admits 
a support of cardinality $<\mathfrak{m}$.
By theorem \ref{net}, for all $X\subseteq A$, $\Sg^{\A}X=\Nr_{\alpha}\Sg^{\B}X$. Seeking a contradiction, we 
assume that such an interpolant does not exist.
Then there exists no interpolant in
$\Sg^{\B}(X_1\cap X_2)$. Indeed, let $c$ be an interpolant in $\Sg^{\B}(X_1\cap X_2)$.
Let $\Gamma=(\beta\sim \alpha)\cap \Delta c$. 
Let $c'={\sf c}_{(\Gamma)}c$. Then $a\leq c'$. Also $b={\sf c}_{(\Gamma)}b$, so that $c'\leq b$.
hence $a\leq c'\leq b$. But $$c'\in \Nr_{\alpha}\Sg^{\B}(X_1\cap X_2)=\Sg^{\Nr_{\alpha}\B}(X_1\cap X_2)
=\Sg^{\A}(X_1\cap X_2).$$
Let $\A_1=\Sg^{\B}X_1$ and $\A_2=\Sg^{\B}X_2$. 
Let $Z_1=\{(J, p): J\subseteq \mathfrak{n}, |J|<\mathfrak{m}, p\in A_1\}$ and define $Z_2$ similarly with $A_2$ replacing $\A_1$.
Then $|Z_1|=|Z_2|\leq \mathfrak{n}$.
To show that $|Z_1|\leq \mathfrak{n}$, let $K$ be a subset of $\mathfrak{n}$ of cardinality $\mathfrak{e}$, the effective degree of $\A_1$.
Then every element $p$ of $\A_1$ is of the form ${\sf s}_{\sigma}q$ with $q\in A_{1K}$ and $\sigma\in {}^{\mathfrak{n}}\mathfrak{n}$.
The number of subsets $J$ of $\mathfrak{n}$ such that $|J|<\mathfrak{m}$ is at most $\sum_{s<\mathfrak{m}}\mathfrak{n}^s=\mathfrak{n}$. 
Let $q\in A_{1K}$ have a support of cardinality $s<\mathfrak{m}$. Then the number of distinct 
elements 
${\sf s}_{\sigma}q$ with $\sigma\in {}^{\mathfrak{n}}\mathfrak{n}$ is at most $\mathfrak{n}^s\leq \mathfrak{n}$. 
Hence there is at most $\mathfrak{n}\cdot\mathfrak{c}$ elements $s_{\sigma}q$ with $\sigma\in {}^{\mathfrak{n}}\mathfrak{n}$ and $q\in A_{1K}$. 
Hence $|Z_1|\leq \mathfrak{n}\cdot\mathfrak{n}\cdot\mathfrak{c}=\mathfrak{n}.$
Let $$\langle (k_i,x_i): i\in \mathfrak{n}\rangle\text {  and  }\langle (l_i,y_i):i\in \mathfrak{n}\rangle$$
be enumerations of $Z_1$ and $Z_2$ respectively, possibly with repititions.
Now there are two functions $u$ and $v$ such that for each $i<\mathfrak{n}$ $u_i, v_i$ are elements of $^\mathfrak{n}\mathfrak{n}$
with 
$$u_i\upharpoonright \mathfrak{n} \sim k_i=Id$$
$$u_i\upharpoonright k_i \text { is one to one }$$
$$x_j\text { and $y_j$ and $a$ and $c$ are independent of $u_i(k_i)$ for all $j\leq i$}$$
$${\sf s}_{u_j}x_j\text {  is independent of $u_i(k_i)$ for all $j<i$}$$
$$v_i\upharpoonright \mathfrak{n}\sim l_i=Id$$
$$v_i\upharpoonright l_i \text { is one to one }$$
$$x_j\text { and $y_j$ and $a$ and $c$ are independent of $v_i(l_i)$ for all $j\leq i$}$$
$${\sf s}_{v_j}y_j\text {  is independent of $v_i(l_i)$ for all $j<i$}$$
and 
$$v_i(l_i)\cap u_j(k_j)=\emptyset\text {  for all  }j\leq i.$$
The existence of such $u$ and $v$ can be proved by by transfinite recursion \cite{MLQ}. Let
$$Y_1= \{a\}\cup \{-{\sf  c}_{(k_i)}x_i\oplus{\sf s}_{u_i}x_i: i\in \mathfrak{n}\},$$
$$Y_2=\{-b\}\cup \{-{\sf  c}_{(l_i)}y_i\oplus {\sf s}_{v_i}y_i:i\in \mathfrak{n}\},$$
$$H_1= fl^{\Rd_{MV}\Sg^{\B}(X_1)}Y_1,\  H_2=fl^{\Rd_{MV}\Sg^\B(X_2)}Y_2,$$ and
$$H=fl^{\Rd_{MV}\Sg^{\B}(X_1\cap X_2)}[(H_1\cap \Sg^{\B}(X_1\cap X_2)
\cup (H_2\cap \Sg^{\B}(X_1\cap X_2)].$$
Then $H$ is a proper filter of $\Sg^{\B}(X_1\cap X_2).$
To prove this it is sufficient to consider any pair of finite, strictly
increasing sequences of ordinals
$$\eta(0)<\eta(1)\cdots <\eta(n-1)<\mathfrak{n}\text { and } \xi(0)<\xi(1)<\cdots
<\xi(m-1)<\mathfrak{n},$$
and to prove that the following condition holds:
\begin{equation}\label{x0-9}
\begin{split}
\textrm{ For any} ~~&b_0, b_1\in \Sg^{\B}(X_1\cap X_2) ~~\textrm{such
that} \\
&a\odot \prod_{i<n}[a^{l-1}\odot (-{\sf  c}_{(k_{\eta(i)})}x_{\eta(i)}\oplus {\sf
s}_{u_{\eta(i)}}x_{\eta(i)})^{l_i}]\leq b_0 \\
\textrm{and}\\
&(-b)\odot \prod_{i<m}[(-b)^{k-1}\odot (-{\sf c}_{(l_{\xi(i)})}y_{\xi(i))}\oplus{\sf
s}_{v_{\xi(i)}}y_{\xi(i)})^{k_i}]\leq b_1\\
\textrm{we have}\\
& b_0\odot b_1\neq 0.
\end{split}
\end{equation}
By induction on $n+m+l-1+k-1=n+m+l+k-2\geq 0$. This can be done by the above argument together with those in \cite{MLQ}.

\item We proceed exactly as above. Proving that $H$ is a proper filter of $\Sg^{\B}(X_1\cap X_2)$,
let $H^*$ be a (proper $MV$) maximal filter of $\Sg^{\B}(X_1\cap X_2)$
containing $H.$
We obtain  maximal filters $F_1$ and $F_2$ of $\Sg^{\B}(X_1)$ and $\Sg^{\B}(X_2)$,
respectively, such that
$$H^*\subseteq F_1,\ \  H^*\subseteq F_2$$
and (*)
$$F_1\cap \Sg^{\B}(X_1\cap X_2)= H^*= F_2\cap \Sg^{\B}(X_1\cap X_2).$$
Now for all $x\in \Sg^{\B}(X_1\cap X_2)$ we have
$$x\in F_1\text { if and only if } x\in F_2.$$
Let $i\in \{1,2\}$. Then $F_i$ by construction satisfies the following:
for each $p\in \B$ and each subset $J\subseteq \mathfrak{n}$ with $|J|<\mathfrak{m}$,  there exists $\rho\in {}^{\mathfrak{n}}\mathfrak{n}$ such that
$$\rho\upharpoonright \mathfrak{n}\sim J=Id_{\mathfrak{n}-J}$$
and $$-{\sf c}_{(J)}p \oplus {\sf s}_{\rho}p\in F_i.$$
Since $${\sf s}_{\rho}p\leq {\sf c}_{(J)}p,$$
we have (**)
$${\sf c}_{(J)}p\in F_i\Longleftrightarrow{\sf s}_{\rho}p\in F_i.$$
Let $\D_i=\Sg^{\A_i}X_i$, $i=1, 2$.
Let $$\psi_i:\D_i\to \F(^{\alpha}\mathfrak{n}, \D_i/F_i)$$ be defined as follows:
$$\psi_i(a)(x)=s_{\bar{x}}^{\B}a/F$$
Note that  $\bar{\tau}=\tau\cup Id_{\beta\sim \alpha}$ is in $^{\mathfrak{n}}\mathfrak{n}$, so that substitutions are evaluated in the big algebra $\B$.
Then, we claim that  $\psi_i$ is a homomorphism. Then using freeness we paste the two maps $\psi_1$ $\psi_2$, obtaining that $\psi=\psi_1\cup \psi_2$ 
is homomorphism from the free algebra to 
$\F(^{\alpha}\mathfrak{n}, [0,1])$ such that $\psi(a-b)\neq 0$ which is a contradiction. As usual,
we check cylindrifiers, and abusing notation for a while, we omit superscripts, in particular, we write $\psi$ instead of $\psi_1$.
Let $x\in {}^{\alpha}\mathfrak{n}$, 
$M\subseteq \alpha,$ $p\in D$.  Then 
$$\psi({\sf c}_{(M)}p)(x)= {\sf s}_x{\sf c}_{(M)}p/F.$$ 
Let $K$ be a support of $p$ such that  $|K|<\mathfrak{m}$
and let $J=M\cap K$. Then
$${\sf c}_{(J)}p={\sf c}_{(M)}p.$$
Let $$\sigma\in {}^\mathfrak{n}\mathfrak{n}, \sigma\upharpoonright \mathfrak{n}\sim J=\bar{\tau}\upharpoonright \mathfrak{n}\sim J,$$
$$\sigma J\cap \tau(K\sim J)=\emptyset$$
and
$$\sigma\upharpoonright J\text { is one to one }.$$
Then $${\sf s}_{x}{\sf c}_{(J)}p={\sf c}_{(\sigma J)}{\sf s}_{\sigma}p.$$
By (**), let $\rho$ be such that 
$$\rho\upharpoonright \sigma J=Id_{\mathfrak{n}\sim \sigma J}$$
and
$${\sf c}_{(\sigma J)}{\sf s}_{\sigma}p\in F\Longleftrightarrow {\sf s}_{\rho}{\sf s}_{\sigma}p\in F
\Longleftrightarrow {\sf s}_{\rho\circ \sigma}p\in F.$$
It follows that 
$${\sf c}_{(\sigma J)}{\sf s}_{\sigma}p/F={\sf s}_{\rho}{\sf s}_{\sigma}p/F={\sf s}_{\rho\circ \sigma}/F.$$
Let $y\in {}^{\mathfrak{n}}\mathfrak{n}$ such that
$$y\upharpoonright M=\rho\circ \sigma \upharpoonright M$$
and (***)
$$y\upharpoonright \mathfrak{n}\sim M=\bar{\tau}\upharpoonright \mathfrak{n}\sim M.$$
If $k\in K\sim M$, $\rho\circ \sigma(k)=\rho\circ \tau(k)$ and since $\tau k\notin \sigma J$ we have
$$\rho \tau(k)=\tau(k)=y(k).$$ So
$$y\upharpoonright K=\rho \circ \sigma \upharpoonright K.$$
Now we have using (**) and (***):
\begin{equation*}
\begin{split}
&\psi({\sf c}_{M}p)x\\
&={\sf s}_x{\sf c}_{(M)}p/F\\
&={\sf s}_x{\sf c}_{(J)}p/F\\
&={\sf c}_{(\sigma J)}{\sf s}_{\sigma}p/F\\
&={\sf s}_{\rho}{\sf s}_{\sigma}p/F\\
&={\sf s}_{\rho\circ \sigma}p/F\\
&={\sf s}_yp/F\\
&\leq {\sf c}_{M}(\psi(p)x\\
\end{split}
\end{equation*}
\end{enumarab}
The other inlusion is left to the reader.
The proof is complete.
\end{demo}

\begin{corollary} Let $\alpha$ be infinite, $G={}^{\alpha}\alpha$ and $T=\wp(\alpha)$. Then every 
$\A\in M_{G,T}$ is representable.
\end{corollary} 
\begin{demo}{Proof} \cite{DM}. Let $\A\in MV_{G,T}$. Let $a$ be non-zero in $\A$. Let $\mathfrak{n}$ and $\B$ an algebra in $\mathfrak{n}$ dimensions such that 
$\A=\Nr_{\alpha}\B.$ Let $Z$ be as in $Z_1$ with $p$ restricted to $A$. Let $((k_i,x_i): i<\mathfrak{n})$ be an enumeration of $Z$. 
Define by transfinite recursion $u$ such that:
$$u_i\upharpoonright \mathfrak{n} \sim k_i=Id$$
$$u_i\upharpoonright k_i \text { is one to one }$$
$$x_j\text { $a$ is independent of $u_i(k_i)$ for all $j\leq i$}$$
$${\sf s}_{u_j}x_j\text {  is independent of $u_i(k_i)$ for all $j<i$}$$
Then form $Y$ as $Y_1$,  take the maximal filter containing $Y$ and $a$, and then define the representation function like $\psi_i$ in the previous proof.
\end{demo} 
\subsection{Interpolation for Pavleka Predicate calculus}

In this section we follow closely \cite{Pav}, however our proof, on the one hand, is much simpler, 
and on the other, it is much more general in at least two respects. We do not resort to the complex constructions of compressions
and dilations resorting to the notions of constants, adressed in the forementioned paper, 
though we use dilations in a disguised form of neat embeddings. 
Furthermore we prove the stronger result of interpolation from which 
we infer the completeness (representability result).
Finally we admit predicate of infinite arity in the corresponding logics.

The Pavelka propositional calculas was introduced to incorporate in the syntax 
truth constants $\bar{r}$ for any $r\in [0,1]$. 
The language then becomes uncountable, but Hajek simplifed this by eliminating from the syntax the irrational truth values. The predicate calculas 
for Pavelka logic was studied by 
Novak, and algebraically by Dragulici and Georgescu. We pursue the latter approach, but in a more general setting.
Let us denote he $MV$ algebra $[0,1]\cap \mathbb{Q}$ by $L$. The following two definitions and lemma are taken from \cite{Pav}.
\begin{definition} A Pavelka algebra is a structure $(\A, \{\bar{r}: r\in L\})$ where $\A$ is an $MV$ algebra and $\{\bar{r}:r\in L\}$ is a subset of $A$ such that 
$\bar{0}=0$, $\bar{r}\oplus \bar{s}=\bar{r\oplus s}$ and $\neg \bar{r}=\bar{\neg r}$, for all $r,s\in L$.  
\end{definition}

\begin{lemma} Let $\B=(\A,\{\bar{r}: r\in L\})$ be a Pavelka algebra. Assume that $P$ is a proper filter of $\B$ and that $r,s\in L$. Then $\bar{r}\in P$ iff 
$r=1$ and $\bar{r}/P\leq \bar{s}/P$ 
iff $r\leq s.$
\end{lemma}

\begin{definition} An existential quantifier on a Pavleka algebra $(\A, \{\bar{r}: r\in L\})$
is an existential quantifier on $\A$ such that $\exists r=r$ for every $r\in L.$
\end{definition}

We consider only substitutions indexed by replacements

\begin{definition} Let $\alpha$ be an ordinal. An $MV$ substitution algebra of dimension $\alpha$, an $SA_{\alpha}$ for short, is an algebra of the form
$${\B}=(\A,\{\bar{r}: r\in L\},{\sf c}_i, {\mathsf s}_i^j)_{i,j< \alpha}$$
where $(\A, \{\bar{r}:r\in L\})$ is a Pavelka algebra and ${\mathsf c}_i, {\mathsf s}_i^j$ are unary operations on $\A$ ($i, j<\alpha$)
such that the ${\sf c}_i$'s are quantifiers and the ${\sf s}_i^j$'s are substitutions. The substitutions
satisfy the same equations as in theorem \ref{d} above.
\end{definition}
For the time being, we restrict our attention to $MV$ substitution algebras which are dimension complemented.
That is any such algebra $\A$  satisfies 
$\alpha\sim \Delta x$ is infinite for all $x\in A.$ In this case $MV$ algebras 
can be viewed as $M_{G,T}$'s with $G$ the semigroup of finite transformations and $T=\wp_{\omega}(\alpha)$ enriched with constant symbols.
The algebra $\Fr_{\beta}^{\rho}SA_{\alpha}$ is defined as before, but we take into consideration the constants when forming subalgebras, 
and we study rather
the pair $(\Fr_{\beta}^{\rho}SA_{\alpha}, \{\bar{r}: r\in L\})$. For $H\subseteq A$ and $a\in A$, set $$[a]_H=\bigvee\{r\in L:\bar{r}\to a\in H\}.$$ 
To prove that  this algebra has the interpolation property,
we proceed as in  the proof of theorem \ref{in}. 

Let $\A$ be the given algebra and $a\in \Sg^{\A}X_1$, $b\in \Sg^{\A}X_2$ such that $a\leq b$. 
By resorting to theorem \ref{net}(1), neatly embed $\A$ into $\Nr_{\alpha}\B$ where $\B$ is of dimension 
$\kappa$, $\kappa$ a regular cardinal $>max(|\alpha|, |A|)$, so that $\A=\Nr_{\alpha}\B$ and 
for $X\subseteq A$, one has $\Sg^{\A}X=\Nr_{\alpha}\Sg^{\B}X.$
If an interpolant exists in $\B$, then an interpolant exists in $\A$, this can be used as the base of the induction to construct the proper filter 
$H$ and from it
maximal filters $F_1$ and $F_2$ in $\Sg^{\B}X_1$ and $\Sg^{\B}X_2,$ respectively.  However, the homomorphisms are defined differently.
We set $\psi_1(b)(x)=[{\sf s}_xp]_{F_1}$ for $b\in \Sg^{\B}X_1$ and $\psi_2$ is defined analogously.  
(Here $x$ is a finite transformation, and ${\sf s}_x$ is the substitution operator 
which is term definable because our algebras are dimension complemented).
Call a maximal filter Henkin, if whenever ${\sf c}_kx\in F$, then ${\sf s}_j^kx\in F$ for some $j\notin \Delta x$.
Both $F_1$ and $F_2$ are Henkin. We need to verify two things. First:
\begin{theorem}\begin{enumarab}
\item  Let $\B\in SA_{\alpha}$. Let $P$ be a Henkin maximal filter. Then the map $\psi:\B\to \F(V, [0,1])$ defined by 
$\psi(p)x= [{\sf s}_{x}p]_P$ is a homomorphism.
\item Let $\A$, $\B$, $F_1, F_2$ be as in the above (sketch of) proof. 
If $a\in \Sg^{\B}(X_1\cap X_2)$, then $[a]_{F_1}=[a]_{F_2}$.
\end{enumarab}
\end{theorem}
\begin{demo}{Proof} We closely follow \cite{Pav}.
\begin{enumarab}
\item  For every $a,b\in A$ the following equalities hold:
\begin{enumerate}
\item $[a]_P=\bigwedge\{r\in L:a\to \bar{r}\in P\}$

\item $[a\oplus b]_P=[a]_P\oplus [b]_P$, $[a\odot b]_P=[a]_P\odot [b]_P$ and $[\neg a]_P=\neg [a]_P$ 
\end{enumerate}
Then it can be easily checked that $$\psi(p\oplus q)(x)=\psi(p)\oplus \psi(a)(x),$$
$$\psi(p)\odot q)(x)=\psi(p)\odot \psi(a)(x)$$
$$\psi (\neg p)(x)=\neg \psi(p)(x).$$
Also $$\psi(\bar{r})(x)=r.$$
We only check cylindrifiers. We have:
$$\psi({\sf c}_ip)(x)=[{\sf s}_x{\sf c}_ip]_P$$
$${\sf c}_i\psi(p)(x)=\bigvee \{\psi(p)(y) : y\equiv_i x\}.$$
We need to show that 
$$[{\sf s}_x{\sf c}_ip]=\bigvee \{{\sf s}_y(p): y\equiv_ix\}$$
We have $p\leq {\sf c}_ip$ so ${\sf s}_yp\leq {\sf s}_x{\sf c}_ip$, so for every $r\in L$, we have $r\to {\sf s}_yp\leq r\to {\sf s}_x{\sf c}_ip$ 
so for every $r\in L$, $r\to {\sf s}_yp\in P$ implies 
$r\to {\sf s}_x{\sf c}_ip$ and so $[{\sf s}_yp]_P\leq [{\sf s}_x{\sf c}_ip]_P$, so $\bigvee {\sf s}_yp\leq {\sf s}_x{\sf c}_ip.$

Asume that for all $y$, $y\equiv_i x$,  we have ${\sf s}_yp<r< {\sf s}_x{\sf c}_ip$. Then
$\exists s>r$ such that $s\to {\sf s}_x{\sf c}_ip$. 
But $r\to {\sf s}_x{\sf c}_ip\geq s\to {\sf s}_x{\sf c}_ip$, it follows that $$r\to {\sf s}_x{\sf c}_ip\in P.$$
This implies that $r\to {\sf s}_{x(i\to u)}\in P.$
This is a contradiction.

\item We proceed as follows. Let $a\in \Sg^{\B}(X_1\cap X_2)$.
Then 
\begin{equation*}
\begin{split}
&a\to \bar{r}\in F_1\\
&\Longleftrightarrow a\to \bar{r}\in F_1\cap \Sg^{\B}(X_1\cap X_2)\\
&\Longleftrightarrow a\to \bar{r}\in F_2\cap \Sg^{\B}(X_1\cap X_2)\\
&\Longleftrightarrow a\to \bar{r}\in F_2\\
\end{split}
\end{equation*}
\end{enumarab}
Then we can deduce that $[a]_{F_1}=[a]_{F_2}$ and we are done.
\end{demo}

 Instead of taking $G=\{[i|j]: i\in \alpha\}$, consisting only of replacements as we did, 
one can take $G$ to be a strongly rich semigroup (with the algebras in question countable) or $G={}^II$ (with
the algebras in question of any cardinality) and in both cases one takes $T=\wp_{\omega}(\alpha)$, then the above proof works verbatim for both cases.
From this we get a substantial generalization of the result in \cite{Pav}, namely

\begin{corollary} If $\A$ is a Pavelka polyadic algebra in the above similarity type, then it is representable, hence the corresponding logic is complete.
\end{corollary}

\section*{Concluding remarks}
\begin{enumarab}
\item  When the algebras considered are Boolean algebras, then we get  the results in \cite{K}, \cite{DM}, \cite{D}, \cite{AUamal}, \cite{MLQ},  \cite{J},  \cite {IGPL}.
 \item Interpolation theorems proved herein, is a typical instance of a Henkin construction \cite{Henkin} formulated algebraically. 
It also has affinity with the works of 
Rasiowa and Sikorski, as far as representability results are concerned.
\item The corollary on representability of $MV$ polyadic algebra is a substantial generalization of the result of Schwartz \cite{Sh}, since we omit the 
condition of local finitenes. 

\item The logic used is a generalization of Kieslers logic \cite{K}, it is its many valued version.

\item The result of Pavelka logic is an extension of the results in \cite{Pav}. It is strictly stronger since predicates are allowed to be infinite.
\item When $G=\wp(\alpha)$ and $T\in \{\wp_{\omega}(\alpha), \wp(\alpha)\}$, then every $MV_{G,T}$ is representable, 
hence the class of representable algebras form a finitely
axiomatizable variety. This variety also has the superamalgamation property. This follows from a well known result of Maksimova \cite{Mak}, 
that interpolation in free algebras of a variety implies 
that the variety has the superamalgmation property. The same applies to the case when $\alpha$ is countable
and $G$ is a strongly rich semigroup.
\item When $G$ is semigroup of finite transformations on $\alpha$, $T=\wp_{\omega}(\alpha)$, then the class of $MV_{G,T}$ algebras 
satisfying $\alpha\sim \Delta x$
in not a variety. It is not hard to see that this class is not closed under products and ultraproducts.  The fact that dimension 
restricted free algebras have the interpolation property can be used in a fairly 
straightforward manner to show that, for $G$ the semigroup of finite transformations on $\alpha+\omega$ and $T=\wp_{\omega}(\alpha+\omega)$, 
the class $\{\A\in MV_{G,T}: \A=\Sg^{\A}\Nr_{\alpha}\A\}$ has the superamalgamation property. If the condition $\alpha\sim \Delta x$ is 
removed then the resulting class is a variety, and it can be shown that it is the same as the class 
${\bf S}\Nr_{\alpha}MV_{\bar{G},\bar{T}}$ where $\bar{G}$ is the semigroup of finite transformation on $\alpha+\omega$ and $T=\wp(\alpha+1).$
We do not know whether this variety is finitely axiomatizable, but it seems unlikely that it is. 

\item Transformation systems can be applied to other algebraisations of fuzzy logic, like $BL$ algebras which are an abstraction of algebras 
arising from $t$ norms.

\item There are two approaches to fuzzy logic. 
The first one is very closely linked with multi-valued logic tradition (Hajek school \cite{H}). 
So a set of designed values is fixed and this enables us to define an entailement relation. 
The deduction apparatus is defined by a suitable set of logical axioms and suitable inference rules. Another approach due to Pavelka and others 
is devoted to defining a deduction apparatus in which approximate reasonings are admitted. 
Such an apparatus is defined by a suitable fuzzy subset of logical axioms and by a suitable set of fuzzy inference rules. 
In the first case the logical consequence 
operator gives the set of logical consequence of a given set of axioms. 
In the latter the logical consequence operator gives the fuzzy subset of logical consequence of a given fuzzy subset of hypotheses.
We saw that both approaches can benefit by applying the methodology of algebraic logic.

\end{enumarab}

\end{document}